\documentclass[11pt]{amsart}
\usepackage[mathscr]{eucal}
\usepackage{amsfonts}
\usepackage{amsmath}
\usepackage{amsthm}
\usepackage{amssymb}
\usepackage{amscd}
\usepackage{latexsym}
\usepackage{slashbox}
\usepackage{lscape}

\usepackage[dvips]{graphicx}
\usepackage{color}

\theoremstyle{plain}
\newtheorem{definition}[equation]{Definition}

\newtheorem{corollary}[equation]{Corollary}
\newtheorem{lemma}[equation]{Lemma}
\newtheorem{proposition}[equation]{Proposition}
\newtheorem{theorem}[equation]{Theorem}

\theoremstyle{definition}
\newtheorem{remark}[equation]{Remark}

\numberwithin{equation}{subsection}
\renewcommand{\mathfrak}{\mathcal}

\begin{document}
\title{The product in the Hochschild cohomology ring of preprojective algebras of Dynkin quivers }
\author{Ching-Hwa Eu}
\address{Department of Mathematics, Massachusetts Institute of
Technology, Cambridge, MA 02139, U.S.A.}
\email{ceu@math.mit.edu}
\maketitle
\pagestyle{myheadings}
\markboth{Ching-Hwa Eu}{The product in the Hochschild cohomology ring of preprojective algebras of Dynkin quivers}
\tableofcontents
\section{Introduction}
In this paper, we compute the product structure of the Hochschild cohomology of preprojective algebras of quivers of type D and E over a field of characteristic zero. This is a continuation of \cite{EE2} where the cohomology spaces together with the grading induced by the natural grading (all arrows have degree $1$) were computed.

Together with the results in \cite{ES2} where it was done for type $A$ (over a  field of any characteristic), this yields a complete description of the product in the Hochschild cohomology ring of preprojective algebras of $ADE$ quivers over a field of characteristic zero.

We note that this description is essentially uniform (i.e. does not refer to particular Dynkin-diagrams), while the proof uses case-by-case arguments.

For our computation, the same complex as in \cite{ES2} is used, namely the one which we get by applying the $Hom$-functor to the Schofield resolution (which is periodic with period $6$) of the algebra.

To compute the cup product, we use the same method as in \cite{ES2}: via the isomorphism $HH^i(A)\equiv \underline{Hom}(\Omega^iA,A)$ (where for an $A$-bimodule $M$ we write $\Omega M$ for the kernel of its projective cover) we identify elements in $HH^i(A)$ with equivalence classes of maps $\Omega^i(A)\rightarrow A$. For $[f]\in HH^i(A)$ and $[g]\in HH^j(A)$, the product is $[f][g]:=[f\circ\Omega^ig]$ in $HH^{i+j}(A)$. We compute all products $HH^i(A)\times HH^j(A)\rightarrow  HH^{i+j}(A)$ for $0\leq i\leq j\leq 5$. The remaining ones follow from the perodicity of the Schofield resolution and the graded commutativity of the multiplication. Some computations are similar to those in \cite{ES2} for type $A$.

In the first part of the paper, we introduce a basis for each cohomology space explicitly (for each quiver). Then in the second part we compute the product in these bases. We use the results about the grading of the cohomology spaces from \cite{EE2} to find the bases and the products.

Note that for connected non-Dynkin quivers, the Hochschild cohomology and its product structure were already calculated in \cite{CBEG} where the situation is much easier because the homological dimension of the preprojective algebra is $2$.

{\bf{Acknowledgements.}} C. Eu wants to thank his advisor P. Etingof and T. Schedler for useful discussions. This work is partially supported by the NSF grant DMS-0504847.

\section{Preliminaries}

\subsection{Quivers and path algebras}
Let $Q$ be a quiver of ADE type with vertex set $I$ and $|I|=r$. 
We write $a\in Q$ to say that $a$ is an arrow in $Q$. 

We define $Q^*$ to be the quiver obtained from $Q$ by reversing
all of its arrows. We call $\bar Q=Q\cup Q^*$ the \emph{double}
of $Q$. 

Let $C$ be the adjacency matrix corresponding to the
quiver $\bar Q$. 

The concatenation of arrows generate the \emph{nontrivial
paths} inside the quiver $\bar Q$. We define $e_i$, $i\in I$ to
be the \emph{trivial path} which starts and ends at $i$. The
\emph{path algebra} $P_{\bar Q}=\mathbb{C}\bar Q$ of $\bar Q$ over $\mathbb{C}$ is the $\mathbb{C}$-algebra with basis the paths in $\bar Q$ and the product $xy$ of two paths $x$ and $y$ to be their concatenation if they are compatible and $0$ if not. We define the \emph{Lie bracket} $[x,y]=xy-yx$.

Let $R=\oplus_{i\in I}\mathbb{C}e_i$. Then $R$ is a commutative semisimple
algebra, and $P_Q$ is naturally an
$R$-bimodule.

\subsection{The preprojective algebra}
Given a quiver $Q$, we define the \emph{preprojective
algebra} $\Pi_Q$ to be the quotient of the path algebra $P_{\bar Q}$ by
the relation $\sum\limits_{a\in Q}[a,a^*]=0$. 

Given a path $x$, we write $x^*$ for the path obtained from $x$ by reversing all arrows.

From now on, we write $A=\Pi_{Q}$.

\subsection{Graded spaces and Hilbert series}

Let $M=\oplus_{d\geq0}M(d)$ be a $\mathbb Z_+$-graded
vector space, with finite dimensional homogeneous subspaces. 
We denote by $M[n]$ the same space with
grading shifted by $n$. The graded dual space $M^*$ is defined by the
formula $M^*(n)=M(-n)^*$.  

\begin{definition} \textnormal{(The Hilbert series of vector spaces)}\\
We define the \emph{Hilbert series} $h_M(t)$ to be the series
\begin{displaymath}
h_M(t)=\sum\limits_{d=0}^{\infty}\dim M(d)t^d.
\end{displaymath}
\end{definition}

\begin{definition} \textnormal{(The Hilbert series of bimodules)}\\
Let $M=\oplus_{d\geq0}M(d)$ be a $\mathbb{Z_+}$-graded bimodule
over the ring $R$, 
so we can write $M=\oplus M_{i,j}$. We define the 
\emph{Hilbert series} $H_M(t)$ to be a matrix valued series with the entries 
\begin{displaymath}
H_M(t)_{i,j}=\sum\limits_{d=0}^{\infty}\dim\ M(d)_{i,j}t^d.
\end{displaymath}
\end{definition}

\subsection{Root system parameters}

Let $w_0$ be the longest element of the Weyl group $W$ of $Q$. 
Then we define $\nu$ to be the involution of $I$, such
that $w_0(\alpha_i)=-\alpha_{\nu(i)}$ (where $\alpha_i$ is the
simple root corresponding to $i\in I$). It turns out that
$\eta(e_i)=e_{\nu(i)}$ (\cite{S}; see \cite{ES2}).

Let $m_i$, $i=1,...,r$, be the exponents of the root system
attached to $Q$, enumerated in the increasing order. 
Let $h=m_r+1$ be the Coxeter number of $Q$. 

Let $P$ be the permutation matrix corresponding to the involution
$\nu$. Let $r_+=\dim\ker(P-1)$ and $r_-=\dim\ker(P+1)$.
Thus, $r_-$ is half the number of vertices which are not fixed by
$\nu$, and $r_+=r-r_-$.

$A$ is finite dimensional, and the following Hilbert series is known from \cite[Theorem 2.3.]{MOV}:

\begin{equation}\label{Hilbert series matrix}
H_A(t)=(1+Pt^h)(1-Ct+t^2)^{-1}.
\end{equation}

We see that the top degree of $A$ is $h-2$, and for the top degree $A^{top}$ part we get the following decomposition in $1$-dimensional submodules:

\begin{equation}\label{topdegree}
A^{top}=A(h-2)=\bigoplus_{i\in I}e_iA(h-2)e_{\nu(i)}
\end{equation}

\section{Hochschild cohomology}
The Hochschild cohomology spaces of $A$ were computed in \cite{EE2}. We recall the results:
\begin{definition}
We define the spaces 
\begin{eqnarray*}
U&=&\oplus_{d<h-2}HH^0(A)(d)[2],\\
L&=&HH^0(A)(h-2),\\
K&=&HH^2(A)[2],\\
Y&=&HH^6(A)(-h-2).
\end{eqnarray*}
\end{definition}

\begin{theorem}
 \begin{enumerate}
 \item $U$ has the following Hilbert series:
  \begin{equation}
   h_U(t)=\sum\limits_{{i=1 \atop m_i<\frac{h}{2}}}^rt^{2m_i}.
  \end{equation}
 \item We have natural isomorphisms 
 \begin{eqnarray*}
 K\equiv \ker (P+1),\\
 L\equiv \ker (P-1),\\
 \end{eqnarray*}
 and
 \[\dim Y=r_+-r_--\#\{i:m_i=\frac{h}{2}\}.\]
 
  \end{enumerate}
\end{theorem}

\begin{theorem}\label{t1}
 For the Hochschild cohomology spaces, we have the following natural isomorphisms:
\begin{align*}
{HH^0}(A)&=U[-2]\oplus L[h-2],\\
{HH^1}(A)&=U[-2],\\
{HH^2}(A)&=K[-2],\\
{HH^3}(A)&=K^*[-2],\\
{HH^4}(A)&=U^*[-2],\\
{HH^5}(A)&=U^*[-2]\oplus Y^*[-h-2],\\
{HH^6}(A)&=U[-2h-2]\oplus Y[-h-2],
\end{align*}
and ${HH^{6n+i}}(A)={HH^i}(A)[-2nh]\,\forall
i\geq1$.
\end{theorem}

\begin{corollary}\label{center}
The center $Z=HH^0(A)$ of $A$ has Hilbert series
$$
h_Z(t)=\sum\limits_{{i=1\atop m_i<\frac{h}{2}}}^rt^{2m_i-2}+r_+t^{h-2}. 
$$
\end{corollary}

\section{Results about the product in the Hochschild cohomology ring of preprojective algebras for D- and E-quivers}
Let $(U[-2])_+$ be the positive degree part of $U[-2]$ (which lies in non-negative degrees). 

We have a decomposition $HH^0(A)=\mathbb{C}\oplus (U[-2])_+\oplus L[-h-2]$ where we have the natural identification $(U[-2])(0)=\mathbb{C}$.

We give a brief description of the product structure in $HH^*(A)$ which will be  computed in this paper. Since the product $HH^i(A)\times HH^j(A)\rightarrow HH^{i+j}(A)$ is graded-commutative, we can assume $i\leq j$ here.

Let $z_0=1\in\mathbb{C}\subset U[-2]\subset HH^0(A)$ (in lowest degree $0$), \\
$\theta_0$ the corresponding element in $HH^1(A)$ (in lowest degree $0$),\\
$\psi_0$ the dual element of $z_0$ in $U^*[-2]\subset HH^5(A)$ (in highest degree $-4$), i.e. $\psi_0(z_0)=1$, \\
$\zeta_0$ the corresponding element in $U^*[-2]\subset HH^4(A)$ (in highest degree $-4$), that is the dual element of $\theta_0$, $\zeta_0(\theta_0)=1$,\\
$\varphi_0:HH^0(A)\rightarrow HH^6(A)$
the natural quotient map (which induces the natural isomorphism $U[-2]\rightarrow U[-2h-2])$ and\\
$\phi$ the quotient map $L\rightarrow Y$ induced by $\varphi_0$.

\begin{theorem}(The product structure in $HH^*(A)$ for quivers of type $D$ and $E$)\\
 \begin{enumerate}
\item
The multiplication by $\varphi_0(z_0)$ induces the natural isomorphisms\\ $\varphi_i:HH^i(A)\rightarrow HH^{i+6}(A)$ $\forall i\geq 1$ and the natural quotient map $\varphi_0$. Therefore, it is enough to compute products $HH^i(A)\times HH^j(A)\rightarrow HH^{i+j}(A)$ with $0\leq i\leq j\leq 5$.
\item The $HH^0(A)$-action on $HH^i(A)$:
\begin{enumerate}
\item ($(U[-2])_+$-action)\\ The action of $(U[-2])_+$ on $U[-2]\subset HH^1(A)$ corresponds to the multiplication 
\begin{eqnarray*}
(U[-2])_+\times U[-2]&\rightarrow& U[-2], \\
(u,v)&\mapsto& u\cdot v
\end{eqnarray*}
in  $HH^0(A)$, projected on $U[-2]\subset HH^0(A)$.\\
$(U[-2])_+$ acts on $U^*[-2]=HH^4(A)$ and $U^*[-2]\subset HH^5(A)$ the following way:
\begin{eqnarray*}
(U[-2])_+\times U^*[-2]&\rightarrow& U^*[-2],\\
(u,f)&\mapsto&u\circ f,
\end{eqnarray*}
where $(u\circ f)(v)=f(uv)$.\\
$(U[-2])_+$ acts by zero on $L[h-2]\subset HH^0(A)$, $HH^2(A)$, $HH^3(A)$ and $Y^*[-h-2]\subset HH^5(A)$.
\item ($L[h-2]$-action)\\
$L[h-2]$ acts by zero on $HH^i(A)$, $1\leq i\leq 4$, and on $U^*[-2]\subset HH^5(A)$. \\
The $L[h-2]$-action on $HH^5(A)$ restricts to
\begin{eqnarray*}
L[h-2]\times Y^*[-h-2]&\rightarrow&U^*[-2],\\
(a,y)&\mapsto& y(\phi(a))\psi_0.
\end{eqnarray*}

\end{enumerate}
\item (Zero products)\\
All products $HH^i(A)\times HH^j(A)\rightarrow HH^{i+j}$, $1\leq i\leq j\leq 5$, where $i+j\geq 6$ or $i,j$ are both odd are zero except the pairings \[HH^1(A)\times HH^5(A)\rightarrow HH^6(A)\] and \[HH^5(A)\times HH^5(A)\rightarrow HH^{10}(A).\]
\item ($HH^1(A)$-products)
\begin{enumerate}
\item The multiplication \[HH^1(A)\times HH^4(A)=U[-2]\times U^*[-2]\rightarrow HH^5(A)\] is the same one as the restriction of \[HH^0(A)\times HH^5(A)\rightarrow HH^5(A)\] on $U[-2]\times U^*[-2]$.\\
\item The multiplication of the subspace $U[-2]_+\subset HH^1(A)$ with $HH^i(A)$ where $i=2,5$ is zero.\\
\item The multiplication by $\theta_0$ induces a symmetric isomorphism \[\alpha:HH^2(A)=K[-2]\rightarrow K^*[-2]=HH^3(A).\] On $HH^5(A)$, it induces a skew-symmetric isomorphism \[\beta:Y^*[-h-2]\rightarrow Y[-h-2]\subset HH^6(A),\] and acts by zero on $U^*[-2]\subset HH^5(A)$. $\alpha$ and $\beta$ will be given by explicit matrices $M_\alpha$ amd $M_\beta$ later.
\end{enumerate}
\item ($HH^2(A)$-products)
\begin{eqnarray*}
HH^2(A)\times HH^2(A)&\rightarrow& HH^4(A),\\ 
(a,b)&\mapsto& \langle-,-\rangle\zeta_0
\end{eqnarray*} 
is given by $\langle-,-\rangle=\alpha$ where $\alpha$ is regarded as a symmetric bilinear form.\\

$HH^2(A)\times HH^3(A)\rightarrow HH^5(A)$ is the multiplication
\begin{eqnarray*}
K[-2]\times K^*[-2]&\rightarrow& HH^5(A),\\
(a,y)&\mapsto&y(a)\psi_0.
\end{eqnarray*}
\item ($HH^5(A)\times HH^5(A)\rightarrow HH^{10}(A)$)\\
The restriction of this product to  \begin{eqnarray*}
Y^*[-h-2]\times Y^*[-h-2]&\rightarrow& HH^{10}(A),\\
 (a,b)&\mapsto& \Omega(-,-)\varphi_4(\zeta_0)
\end{eqnarray*}  
is given by $\Omega(-,-)=-\beta$ where $\beta$ is regarded as a skew-symmetric bilinear form.

The multiplication of the subspace $U^*[-2]\subset HH^5(A)$ with $HH^5(A)$ is zero.
\end{enumerate}
\end{theorem}
\section{Some basic facts about preprojective algebras}

\subsection{Labeling of quivers}
From now on, we use the following labellings for the different types of quivers:

\subsubsection{$Q=D_{n+1}$}
\begin{figure}[htp]
\input{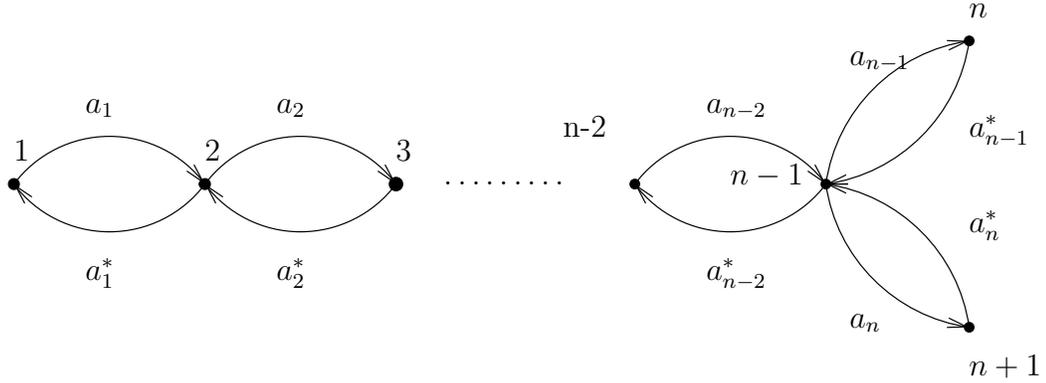}
\caption{$D_{n+1}$-quiver}
\end{figure}
$A$ is the path algebra modulo the relations
\begin{eqnarray*}
a_1^*a_1&=&0,\\
a_{i+1}^*a_{i+1}&=&a_ia_i^*,\quad 1\leq n-3\\
a_{n-1}^*a_{n-1}+a_n^*a_n&=&a_{n-2}a_{n-2}^*.
\end{eqnarray*}

\subsubsection{$Q=E_6$}
\begin{figure}[htp]
\input{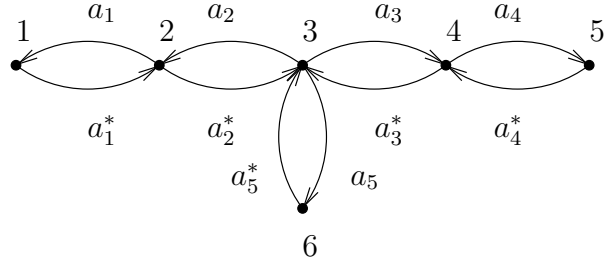}
\caption{$E_6$-quiver}
\end{figure}
$A$ is the path algebra modulo the relations
\begin{eqnarray*}
a_1a_1^*=a_4a_4^*=a_5a_5^*&=&0,\\
a_1^*a_1&=&a_2a_2^*,\\
a_4^*a_4&=&a_3a_3^*,\\
a_2^*a_2+a_3^*a_3+a_5a_5^*&=&0.
\end{eqnarray*}

\subsubsection{$Q=E_7$}
\begin{figure}[htp]
\input{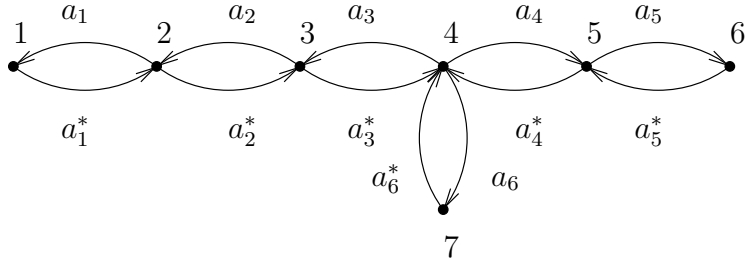}
\caption{$E_7$-quiver}
\end{figure}
$A$ is the path algebra modulo the relations
\begin{eqnarray*}
a_1a_1^*=a_5a_5^*=a_6a_6^*&=&0,\\
a_1^*a_1&=&a_2a_2^*,\\
a_2^*a_2&=&a_3a_3^*,\\
a_5^*a_5&=&a_4a_4^*,\\
a_3^*a_3+a_4^*a_4+a_6a_6^*&=&0.
\end{eqnarray*}

\begin{figure}[htp]
\input{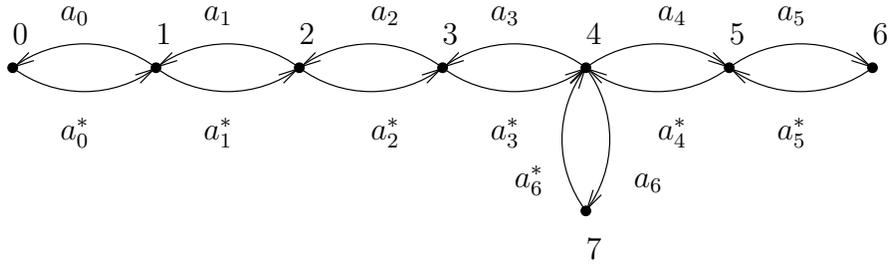}
\caption{$E_8$-quiver}
\end{figure}
$A$ is the path algebra modulo the relations
\begin{eqnarray*}
a_0a_0^*=a_5a_5^*=a_6a_6^*&=&0,\\
a_0^*a_0&=&a_1a_1^*,\\
a_1^*a_1&=&a_2a_2^*,\\
a_2^*a_2&=&a_3a_3^*,\\
a_5^*a_5&=&a_4a_4^*,\\
a_3^*a_3+a_4^*a_4+a_6a_6^*&=&0.
\end{eqnarray*}

\subsection{Preprojective algebras by numbers}
We summarize useful numbers associated to preprojective algebras, by quiver:\\
 \begin{tabular}{|c|c|c|c|c|}
 \hline
 $Q$& exponents $m_i$&$h$&$\deg A^{top}$& degrees $HH^0(A)$\\
 \hline
 $D_{n+1}$ \begin{tabular}{l} $n$ odd\\$n$ even\end{tabular} & $1,3,\ldots, 2n-1,n$  & $2n$ & $2n-2$ & \begin{tabular}{l}
                                                                                     $0,4,\ldots,2n-6,2n-2$\\$0,4,\ldots 2n-4, 2n-2$ 
                                                                                     \end{tabular}\\
\hline
$E_6$ & $1,4,5,7,8,11$ & $12$ & $10$ & $0,6,8,10$\\
\hline
$E_7$ & $1,5,7,9,11,13,17$ & $18$ & $16$ & $0,8,12,16$\\
\hline 
$E_8$ & $1,7,11,13,17,19,23,29$ & $30$ & $28$ & $0,12,20,24,28$\\
\hline
\end{tabular}

We see that for quivers of type $D$ and $E$, the degrees of the space $U$ (which are $2 m_i,\,m_i<\frac{h}{2}$) are even and range from $0$ to $h-2$. 

We get the following degree ranges for the Hochschild cohomology:
\[
\begin{array}{rlc}
{HH^0}(A)&=U[-2]\oplus L[h-2] & 0\leq\deg HH^0(A)\leq h-2\\
{HH^1}(A)&=U[-2] & 0\leq\deg HH^1(A)\leq h-4\\
{HH^2}(A)&=K[-2] & \deg HH^2(A)=-2\\
{HH^3}(A)&=K[-2] & \deg HH^3(A)=-2\\
{HH^4}(A)&=U^*[-2] & -h\leq\deg HH^4(A)\leq -4\\
{HH^5}(A)&=U^*[-2]\oplus Y^*[-h-2] & -h-2\leq \deg HH^5(A)\leq -4\\
{HH^6}(A)&=U[-2h-2]\oplus Y[-h-2] & -2h\leq\deg HH^6(A)\leq -h-2
\end{array}
\]

\subsection{Frobenius algebras and Nakayama automorphism}\label{Frobenius}
\begin{definition}
Let $\mathfrak{A}$ be a finite dimensional unital $\mathbb{C}-$algebra. We
call it Frobenius if there is a linear function $f:\mathfrak{A}\rightarrow\mathbb{C}$, such that the form $(x,y):=f(xy)$ is nondegenerate, or, equivalently, if there exists an isomorphism $\phi:\mathfrak{A}\stackrel{\simeq}{\rightarrow}\mathfrak{A}^*$ of left $\mathfrak{A}-$modules: given $f$, we can define $\phi(a)(b)=f(ba)$, and given $\phi$, we define $f=\phi(1)$.
\end{definition}
\begin{remark}
If $\tilde f$ is another linear function satisfying the same
properties as $f$ from above, then $\tilde f(x)=f(xa)$ for some
invertible $a\in \mathfrak{A}$. Indeed, we define the form $\{a,b\}=\tilde f(ab)$. Then $\{-,1\}\in \mathfrak{A}^*$, so there is an $a\in \mathfrak{A}$, such that $\phi(a)=\{-,1\}$. Then $\tilde f(x)=\{x,1\}=\phi(a)(x)=f(xa)$.
\end{remark}

\begin{definition} 
Given a Frobenius algebra $\mathfrak{A}$ (with a function $f$ inducing a
bilinear form $(-,-)$ from above), the automorphism
$\eta:\mathfrak{A}\rightarrow \mathfrak{A}$ defined by the equation $(x,y)=(y,\eta(x))$
is called the \emph{Nakayama automorphism} (corresponding to
$f$).
\end{definition}

\begin{remark}We note that the freedom in choosing $f$ implies 
that $\eta$ is uniquely determined up to an inner automorphism. Indeed,
let $\tilde f(x)=f(xa)$ and define the bilinear form $\{a,b\}=\tilde f(ab)$. Then
\begin{align*}
\{x,y\}&=\tilde f(xy)=f(xya)=(x,ya)=(ya,\eta(x))=f(ya\eta(x)a^{-1}a)\\
&=(y,a\eta(x)a^{-1}).
\end{align*}
\end{remark}

It is known that $A$ is a Frobenius algebra (see e.g. \cite{ES2},\cite{MOV}). 

The linear function $f:A\rightarrow\mathbb{C}$ is zero in the non-top degree part of $A$. It maps a top degree element $\omega_i\in e_iA^{top}e_{\nu(i)}$ to $1$. 
It is uniquely determined by the choice of one of these $\omega_i$ and a Nakayama automorphism. 

For each quiver, we define a Nakayama automorphism and make a choice of one $\omega_i\in e_iA^{top}e_{\nu(i)}$:

\subsubsection{$Q=D_{n+1}$, $n$ odd}
We define $\eta$ by
\begin{eqnarray}
  \eta(a_{i})&=&-a_{i},\\
  \eta(a_{i}^*)&=&a_i^*,\\
  \end{eqnarray}
and  
  \begin{equation}
  \omega_1=a_1^*\ldots a_{n-2}^*a_{n-1}^*a_{n-1}a_{n-2}\ldots a_1.
\end{equation}
\subsubsection{$Q=D_{n+1}$, $n$ even}
We define $\eta$ by
\begin{eqnarray}
  \forall i\leq n-2:\quad\eta(a_{i})&=&-a_{i},\\
  \forall i\leq n-2:\quad\eta(a_{i}^*)&=&a_i^*,\\
  \eta(a_{n-1})&=&-a_{n},\\
  \eta(a_{n-1}^*)&=&a_n^*,\\
  \eta(a_{n})&=&-a_{n-1},\\
  \eta(a_{n}^*)&=&a_{n-1}^*,
  \end{eqnarray}
  \begin{equation}
  \omega_1=a_1^*\ldots a_{n-2}^*a_{n-1}^*a_{n-1}a_{n-2}\ldots a_1.
\end{equation}

\subsubsection{$Q=E_6$}
We define $\eta$ by
\begin{eqnarray}
 \eta(a_1)&=&-a_4,\\
 \eta(a_1^*)&=&a_4^*,\\
 \eta(a_2)&=&-a_3,\\
 \eta(a_2^*)&=&a_3^*,\\
 \eta(a_5)&=&-a_5,\\
 \eta(a_5^*)&=&a_5^*,
 \end{eqnarray}
 and
 \begin{equation}
 \omega_3=a_3^*a_3(a_2^*a_2a_3^*a_3)^2.
\end{equation}

\subsubsection{$Q=E_7$}
We define $\eta$ by
\begin{eqnarray}
\eta(a_i)&=&-a_i,\\
\eta(a_i^*)&=&a_i^*,
\end{eqnarray}
and
\begin{equation}
\omega_4=(a_4^*a_4a_3^*a_3)^4.
\end{equation}

\subsubsection{$Q=E_8$}
We define $\eta$ by
\begin{eqnarray}
\eta(a_i)&=&-a_i,\\
\eta(a_i^*)&=&a_i^*,
\end{eqnarray}
and
\begin{equation}
\omega_4=(a_4^*a_4a_3^*a_3)^7
\end{equation}

\subsection{The Schofield resolution}
We recall the Schofield resolution of $A$ from \cite{S}.

Define the 
$A-$bimodule $\mathfrak{N}$ obtained from $A$ by twisting the
right action by $\eta$, i.e., 
$\mathfrak{N}=A$ as a vector space, 
and $\forall a,b\in A,x\in\mathfrak{N}:a\cdot x\cdot b=ax\eta(b).$
Introduce the notation $\epsilon_a=1$ if $a\in Q$,
$\epsilon_a=-1$ if $a\in Q^*$. Let $x_i$ be a homogeneous basis
of $A$ and $x_i^*$
the dual basis under the form attached to the Frobenius algebra
$A$. Let $V$ be the bimodule spanned by the edges of $\bar Q$. 

We start with the following exact sequence:
\[
0\rightarrow\mathfrak{N}[h]\stackrel{i}{\rightarrow}P_2\stackrel{d_2}{\rightarrow}P_1\stackrel{d_1}{\rightarrow}P_0\stackrel{d_0}{\rightarrow}A\rightarrow 0,
\]
where $P_2=A\otimes_R A[2]$, $P_1=A\otimes_RV\otimes_RA$, $P_0=A\otimes_RA$, 
\begin{align*}
d_0(x\otimes y)&=xy,\\
d_1(x\otimes v\otimes y)&=xv\otimes y-x\otimes vy,\\
d_2(z\otimes t)&=\sum\limits_{a\in\bar Q}\epsilon_aza\otimes a^*\otimes t+\sum\limits_{a\in\bar Q}\epsilon_az\otimes a\otimes a^*t,\\
i(a)&=a\sum x_i\otimes x_i^*.
\end{align*}

Since $\eta^2=1,$ we can make a canonical identification
$A=\mathfrak{N}\otimes_A\mathfrak{N}$ (via $x\mapsto x\otimes
1$), so by tensoring the above exact sequence with
$\mathfrak{N}$, connecting with the original exact sequence and repeating this process, we get the Schofield resolution
\[
 \ldots\rightarrow P_6\stackrel{d_6}{\rightarrow}P_5\stackrel{d_5}{\rightarrow}P_4\stackrel{d_4}{\rightarrow}P_3\stackrel{d_3}{\rightarrow}P_2\stackrel{d_2}{\rightarrow}P_1\stackrel{d_1}{\rightarrow}P_0\stackrel{d_0}{\rightarrow}A\rightarrow 0,
\]
with
\[P_{i+3}=(P_i\otimes_R\mathfrak{N})[h].\]

We will work with the Hochschild cohomology complex obtained from this resolution, which is given explicitly in \cite[Subsection 4.5.]{EE2}.

\section{Basis and Hilbert series}\label{basis D_{n+1}}

We need to work with the Hilbert series and with an explicit basis of $A$. We do this for each type of quiver separately. 

We write $B$ for a set of all basis elements of $A$, $B_{i,-}$ for a basis of $e_iA$, $B_{-,j}$ for a basis of $Ae_j$,
$B_{i,j}$ for a basis of $e_iAe_j$ and $B_{i,j}(d)$ for a basis of $e_iAe_j(d)$.\\

\subsection{$Q=D_{n+1}$}

A basis of $A$ is given by the following elements:

For \underline{$k,j\leq n-1$}:
\begin{eqnarray*}
 B_{k,n}&=&\{(a_{k-1}a_{k-1}^*)^la_k^*\cdots a_{n-2}^*a_{n-1}^*|0\leq l\leq k-1\},\\
 B_{k,n+1}&=&\{(a_{k-1}a_{k-1}^*)^la_k^*\cdots a_{n-2}^*a_{n}^*|0\leq l\leq k-1\},\\
 B_{n,n}&=&\{(a_{n-1}a_n^*a_na_{n-1}^*)^l|0\leq l\leq\left\{\begin{array}{ll}\frac{n-1}{2}& n\,odd,\\\frac{n-2}{2} & n\,even\end{array}\right.\},\\
 B_{n+1,n+1}&=&\{(a_na_{n-1}^*a_{n-1}a_n^*)^l|0\leq l\leq\left\{\begin{array}{ll} \frac{n-1}{2} & n\,odd,\\\frac{n-2}{2} & n\,even\end{array}\right.\},\\
 B_{n+1,n}&=&\{a_na_{n-1}^*(a_{n-1}a_n^*a_na_{n-1}^*)^l|0\leq l\leq\left\{\begin{array}{ll}\frac{n-3}{2}& n\,odd,\\\frac{n-2}{2} & n\,even\end{array}\right.\},\\
 B_{n,n+1}&=&\{a_{n-1}a_n^*(a_na_{n-1}^*a_{n-1}a_n^*)^l|0\leq l\leq\left\{\begin{array}{ll} \frac{n-3}{2}& n\,odd,\\\frac{n-2}{2} & n\,even\end{array}\right.\},\\
 B_{n,j}&=&\{a_{n-1}a_{n-2}\cdots a_j(a_{j-1}a_{j-1}^*)^l|0\leq l\leq j-1\},\\
 B_{n+1,j}&=&\{a_{n}a_{n-2}\cdots a_j(a_{j-1}a_{j-1}^*)^l|0\leq l\leq j-1\}.
\end{eqnarray*}
For \underline{$k\leq j\leq n-1$},
\begin{eqnarray*}
B_{k,j}&=&\{(a_{k-1}a_{k-1}^*)^la_k^*\cdots a_{j-1}^*|0\leq l\leq\min\{k-1,n-j-1\}\}\cup\\
&&\{(a_{k-1}a_{k-1}^*)^la_k^*\cdots a_{n-1}^*a_{n-1}a_{n-2}a_j|0\leq l\leq k-1\}\cup\\
&&\{(a_{k-1}a_{k-1}^*)^la_k^*\cdots a_{n}^*a_{n}a_{n-2}a_j|0\leq l\leq k-1+j-n\}.
\end{eqnarray*}
For \underline{$j<k\leq n-1$},
\begin{eqnarray*}
B_{k,j}&=&\{a_{k-1}\cdots a_j(a_j^*a_j)^l|0\leq l\leq\min\{n-k-1,j-1\}\}\cup\\
&&\{a_k^*\cdots a_{n-2}^*a_{n-1}^*a_{n-1}a_{n-2}\cdots a_j(a_j^*a_j)^l|0\leq l\leq j-1\}\cup\\
&&\{a_k^*\cdots a_{n-2}^*a_{n}^*a_{n}a_{n-2}\cdots a_j(a_j^*a_j)^l|0\leq l\leq j-1+k-n\}.
\end{eqnarray*}

\subsection{$Q=E_6$}

We give the columns of the Hilbert series $H_A(t)$ which can be calculated from \ref{Hilbert series matrix}:
\[
(H_A(t)_{i,1})_{1\leq i\leq 6}=\left(
\begin{array}{l}
1+t^6\\
t+t^5+t^7\\
t^2+t^4+t^6+t^8\\
t^3+t^5+t^9\\
t^4+t^{10}\\
t^3+t^7
\end{array}
\right),
\]

\[
(H_A(t)_{i,2})_{1\leq i\leq 6}=\left(
\begin{array}{l}
t+t^5+t^7\\
1+t^2+t^4+2t^6+t^8\\
t+2t^3+2t^5+2t^7+t^9\\
t^2+2t^4+t^6+t^8+t^{10}\\
t^3+t^5+t^9\\
t^2+t^4+t^6+t^8
\end{array}
\right),
\]

\[
(H_A(t)_{i,3})_{1\leq i\leq 6}=\left(
\begin{array}{l}
t^2+t^4+t^6+t^8\\
t+2t^3+2t^5+2t^7+t^9\\
1+2t^2+3t^4+3t^6+2t^8+t^{10}\\
t+2t^3+2t^5+2t^7+t^9\\
t^2+t^4+t^6+t^8\\
t+t^3+2t^5+t^7+t^9
\end{array}
\right),
\]

\[
(H_A(t)_{i,4})_{1\leq i\leq 6}=\left(
\begin{array}{l}
t^3+t^5+t^9\\
t^2+2t^4+t^6+t^8+t^{10}\\
t+2t^3+2t^5+2t^7+t^9\\
1+t^2+t^4+2t^6+t^8\\
t+t^5+t^7\\
t^2+t^4+t^6+t^8
\end{array}
\right),
\]

\[
(H_A(t)_{i,5})_{1\leq i\leq 6}=\left(
\begin{array}{l}
t^4+t^{10}\\
t^3+t^5+t^9\\
t^2+t^4+t^6+t^8\\
t+t^5+t^7\\
1+t^6\\
t^3+t^7
\end{array}
\right),
\]

\[
(H_A(t)_{i,6})_{1\leq i\leq 6}=\left(
\begin{array}{l}
t^3+t^7\\
t^2+t^4+t^6+t^8\\
t+t^3+2t^5+t^7+t^9\\
t^2+t^4+t^6+t^8\\
t^3+t^7\\
1+t^4+t^6+t^{10}
\end{array}
\right).
\]

\subsection{$Q=E_7$}

We give the columns of the Hilbert series matrix $H_A(t)$ of $A$ which can be calculated from \ref{Hilbert series matrix}:
 
\[
(H_A(t)_{i,1})_{1\leq i\leq 7}=\left(
\begin{array}{l}
  1+t^8+t^{16}\\
  t+t^7+t^9+t^{15}\\
  t^2+t^6+t^8+t^{10}+t^{14}\\
  t^3+t^5+t^7+t^9+t^{11}+t^{13}\\ 
  t^4+t^6+t^{10}+t^{12}\\
  t^5+t^{11}\\
  t^4+t^8+t^{12}
\end{array}
\right),
\]

\[
(H_A(t)_{i,2})_{1\leq i\leq 7}=\left(
\begin{array}{l}
t+t^7+t^9+t^{15}\\
1+t^2+t^6+2t^8+t^{10}+t^{14}+t^{16}\\ 
t+t^3+t^5+2t^7+2t^9+t^{11}+t^{13}+t^{15}\\
t^2+2t^4+2t^6+2t^8+2t^{10}+2t^{12}+t^{14}\\
t^3+2t^5+t^7+t^9+2t^{11}+t^{13}\\
t^4+t^6+t^{10}+t^{12}\\
t^3+t^5+t^7+t^9+t^{11}+t^{13}
\end{array}
\right),
\]

\[
(H_A(t)_{i,3})_{1\leq i\leq 7}=\left(
\begin{array}{l}
t^2+t^6+t^8+t^{10}+t^{14}\\
t+t^3+t^5+2t^7+2t^9+t^{11}+t^{13}+t^{15}\\
1+t^2+2t^4+3t^6+3t^8+2t^{10}+2t^{12}+t^{14}+t^{16}\\
1+2t^2+3t^4+4t^6+4t^8+4t^{10}+3t^{12}+2t^{14}+t^{16}\\
t+2t^3+2t^5+3t^7+3t^9+2t^{11}+2t^{13}+t^{15}\\
t^2+t^4+t^6+2t^8+t^{10}+t^{12}+t^{14}\\
t+t^3+2t^5+2t^7+2t^9+2t^{11}+t^{13}+t^{15}\\
\end{array}
\right)
\]

\[
(H_A(t)_{i,4})_{1\leq i\leq 7}=\left(
\begin{array}{l}
t^3+t^5+t^7+t^9+t^{11}+t^{13}\\
t^2+2t^4+2t^6+2t^8+2t^{10}+2t^{12}+t^{14}\\
t+2t^3+3t^5+3t^7+3t^9+3t^{11}+2t^{13}+t^{15}\\
1+2t^2+3t^4+4t^6+4t^8+4t^{10}+3t^{12}+2t^{14}+t^{16}\\
t+2t^3+2t^5+3t^7+3t^9+2t^{11}+2t^{13}+t^{15}\\
t^2+t^4+t^6+2t^8+t^{10}+t^{12}+t^{14}\\
t+t^3+2t^5+2t^7+2t^9+2t^{11}+t^{13}+t^{15}
\end{array}
\right),
\]

\[
(H_A(t)_{i,5})_{1\leq i\leq 7}=\left(
\begin{array}{l}
t^4+t^6+t^{10}+t^{12}\\
t^3+2t^5+t^7+t^9+2t^{11}+t^{13}\\
t^2+2t^4+2t^6+2t^8+2t^{10}+2t^{12}+t^{14}\\
t+2t^3+2t^5+3t^7+3t^9+2t^{11}+2t^{13}+t^{15}\\
1+t^2+t^4+2t^6+2t^8+2t^{10}+t^{12}+t^{14}+t^{16}\\
t+t^5+t^7+t^9+t^{11}+t^{15}\\
t^2+t^4+t^6+2t^8+t^{10}+t^{12}+t^{14}
\end{array}
\right),
\]

\[
(H_A(t)_{i,6})_{1\leq i\leq 7}=\left(
\begin{array}{l}
t^5+t^{11}\\
t^4+t^6+t^{10}+t^{12}\\
t^3+t^5+t^7+t^9+t^{11}+t^{13}\\
t^2+t^4+t^6+2t^8+t^{10}+t^{12}+t^{14}\\
t+t^5+t^7+t^9+t^{11}+t^{15}\\
1+t^6+t^{10}+t^{16}\\
t^3+t^7+t^9+t^{13}
\end{array}
\right),
\]

\[
(H_A(t)_{i,7})_{1\leq i\leq 7}=\left(
\begin{array}{l}
t^4+t^8+t^{12}\\
t^3+t^5+t^7+t^9+t^{11}+t^{13}\\
t^2+t^4+2t^6+t^8+2t^{10}+t^{12}+t^{14}\\
t+t^3+2t^5+2t^7+2t^9+2t^{11}+t^{13}+t^{15}\\
t^2+t^4+t^6+2t^8+t^{10}+t^{12}+t^{14}\\
t^3+t^7+t^9+t^{13}\\
1+t^4+t^6+t^8+t^{10}+t^{12}+t^{16}
\end{array}
\right).
\]

\newpage

\subsection{$Q=E_8$}

We give the columns of the Hilbert series matrix $H_A(t)$ of $A$ which can be calculated from \ref{Hilbert series matrix}:

\begin{eqnarray*}\lefteqn{(H_A(t)_{i,1})_{1\leq i\leq8}=}\\
&=\left(
 \begin{array}{l}
 1+t^{10}+t^{18}+t^{28}\\
 t+t^9+t^{11}+t^{17}+t^{19}+t^{27}\\
 t^2+t^8+t^{10}+t^{12}+t^{16}+t^{18}+t^{20}+t^{26}\\
 t^3+t^7+t^9+t^{11}+t^{13}+t^{15}+t^{17}+t^{19}+t^{21}+t^{25}\\
 t^4+t^6+t^8+t^{10}+t^{12}+2t^{14}+t^{16}+t^{18}+t^{20}+t^{22}+t^{24}\\
 t^5+t^7+t^{11}+t^{13}+t^{15}+t^{17}+t^{21}+t^{23}\\
 t^6+t^{12}+t^{16}+t^{22}\\
 t^5+t^9+t^{13}+t^{15}+t^{19}+t^{23}
 \end{array}
\right)
\end{eqnarray*}

\begin{eqnarray*}\lefteqn{(H_A(t)_{i,2})_{1\leq i\leq8}=}\\
&=\left(
 \begin{array}{l}
 t+t^9+t^{11}+t^{17}+t^{19}+t^{27}\\
 1+t^2+t^8+2t^{10}+t^{12}+t^{16}+2t^{18}+t^{20}+t^{26}+t^{28}\\
 t+t^3+t^7+2t^9+2t^{11}+t^{13}+t^{15}+2t^{17}+2t^{19}+t^{21}+t^{25}+t^{27}\\
 t^2+t^4+t^6+2t^8+2t^{10}+2t^{12}+2t^{14}+2t^{16}+2t^{18}+2t^{20}+t^{22}+t^{24}+t^{26}\\
 t^3+2t^5+2t^7+2t^9+2t^{11}+3t^{13}+3t^{15}+2t^{17}+2t^{21}+2t^{23}+t^{25}\\
 t^4+2t^6+t^8+t^{10}+2t^{12}+2t^{14}+2t^{16}+t^{18}+t^{20}+2t^{22}+t^{24}\\
 t^5+t^7+t^{11}+t^{13}+t^{15}+t^{17}+t^{21}+t^{23}\\
 t^4+t^6+t^8+t^{10}+t^{12}+2t^{14}+t^{16}+t^{18}+t^{20}+t^{22}+t^{24}
 \end{array}
\right),
\end{eqnarray*}  

\[(H_A(t)_{i,3})_{1\leq i\leq8}=\left(
 \begin{array}{l}
 t^2+t^8+t^{10}+t^{12}+t^{16}+t^{18}+t^{20}+t^{26}\\
 t+t^3+t^7+2t^9+2t^{11}+t^{13}+t^{15}\\\quad+2t^{17}+2t^{19}+t^{21}+t^{25}+t^{27}\\
 1+t^2+t^4+t^6+2t^8+3t^{10}+2t^{12}+2t^{14}\\\quad+2t^{16}+3t^{18}+2t^{20}+t^{22}+t^{24}+t^{26}+t^{28}\\
 t+t^3+2t^5+2t^7+3t^9+3t^{11}+3t^{13}+3t^{15}\\\quad+3t^{17}+3t^{19}+2t^{21}+2t^{23}+t^{25}+t^{27}\\
 t^2+2t^4+3t^6+3t^8+3t^{10}+4t^{12}+4t^{14}\\\quad+4t^{16}+3t^{18}+3t^{20}+3t^{22}+2t^{24}+t^{26}\\
 t^3+2t^5+2t^7+2t^9+2t^{11}+3t^{13}+3t^{15}\\\quad+2t^{17}+2t^{19}+2t^{21}+2t^{23}+t^{25}\\
 t^4+t^6+t^8+t^{10}+t^{12}+2t^{14}\\\quad+t^{16}+t^{18}+t^{20}+t^{22}+t^{24}\\
 t^3+t^5+2t^7+t^9+2t^{11}+2t^{13}+2t^{15}\\\quad+2t^{17}+t^{19}+2t^{21}+t^{23}+t^{25}
 \end{array}
\right)\]

\[(H_A(t)_{i,4})_{1\leq i\leq8}=\left(
 \begin{array}{l}
t^3+t^7+t^9+t^{11}+t^{13}+t^{15}\\\quad+t^{17}+t^{19}+t^{21}+t^{25}\\
t^2+t^4+t^6+2t^8+2t^{10}+2t^{12}+2t^{14}\\\quad+2t^{16}+2t^{18}+2t^{20}+t^{22}+t^{24}+t^{26}\\
t+t^3+2t^5+2t^7+3t^9+3t^{11}+3t^{13}+3t^{15}\\\quad+3t^{17}+3t^{19}+2t^{21}+2t^{23}+t^{25}+t^{27}\\
1+t^2+2t^4+3t^6+3t^8+4t^{10}+4t^{12}+4t^{14}\\\quad+4t^{16}+4t^{18}+3t^{20}+3t^{22}+2t^{24}+t^{26}+t^{28}\\
t+2t^3+3t^5+4t^7+4t^9+5t^{11}+5t^{13}+5t^{15}\\\quad+5t^{17}+4t^{19}+4t^{21}+3t^{23}+2t^{25}+t^{27}\\
t^2+2t^4+2t^6+3t^8+3t^{10}+3t^{12}+4t^{14}\\\quad+3t^{16}+3t^{18}+3t^{20}+2t^{22}+2t^{24}+t^{26}\\
t^3+t^5+t^7+2t^9+t^{11}+2t^{13}+2t^{15}\\\quad+t^{17}+2t^{19}+t^{21}+t^{23}+t^{25}\\
t^2+t^4+2t^6+2t^8+2t^{10}+3t^{12}+2t^{14}\\\quad+3t^{16}+2t^{18}+2t^{20}+2t^{22}+t^{24}+t^{26} 
 \end{array}
\right)\]

\[(H_A(t)_{i,5})_{1\leq i\leq8}=\left(
 \begin{array}{l}
t^4+t^6+t^8+t^{10}+t^{12}+2t^{14}\\\quad+t^{16}+t^{18}+t^{20}+t^{22}+t^{24}\\
t^3+2t^5+2t^7+2t^9+2t^{11}+3t^{13}+3t^{15}\\\quad+2t^{17}+2t^{19}+2t^{21}+2t^{23}+t^{25}\\
t^2+2t^4+3t^6+3t^8+3t^{10}+4t^{12}+4t^{14}\\\quad+4t^{16}+3t^{18}+3t^{20}+3t^{22}+2t{24}+t^{26}\\
t+2t^3+3t^5+4t^7+4t^9+5t^{11}+5t^{13}+5t^{15}\\\quad+5t^{17}+4t^{19}+4t^{21}+3t^{23}+2t^{25}+t^{27}\\
1+2t^2+3t^4+4t^6+5t^8+6t^{10}+6t^{12}+6t^{14}\\\quad+6t^{16}+6t^{18}+5t^{20}+4t^{22}+3t^{24}+2t^{26}+t^{28}\\
t+2t^3+2t^5+3t^7+4t^9+4t^{11}+4t^{13}+4t^{15}\\\quad+4t^{17}+4t^{19}+3t^{21}+2t^{23}+2t^{25}+t^{27}\\
t^2+t^4+t^6+2t^8+2t^{10}+2t^{12}+2t^{14}\\\quad+2t^{16}+2t^{18}+2t^{20}+t^{22}+t^{24}+t^{26}\\
t+t^3+2t^5+2t^7+3t^9+3t^{11}+3t^{13}+3t^{15}\\\quad+3t^{17}+3t^{19}+2t^{21}+2t^{23}+t^{25}+t^{27}
 \end{array}
\right)\]

\[
(H_A(t)_{i,6})_{1\leq i\leq8}=\left(
 \begin{array}{l}
t^5+t^7+t^{11}+t^{13}+t^{15}+t^{17}+t^{21}+t^{23}\\
t^4+2t^6+t^8+t^{10}+2t^{12}+2t^{14}\\\quad+2t^{16}+t^{18}+t^{20}+2t^{22}+t^{24}\\
t^3+2t^5+2t^7+2t^9+2t^{11}+3t^{13}+3t^{15}\\\quad+2t^{17}+2t^{19}+2t^{21}+2t^{23}+t^{25}\\
t^2+2t^4+2t^6+3t^8+3t^{10}+3t^{12}+4t^{14}+3t^{16}\\\quad+3t^{18}+3t^{20}+2t^{22}+2t^{24}+t^{26}\\
t+2t^3+2t^5+3t^7+4t^9+4t^{11}+4t^{13}\\\quad+4t^{15}+4t^{17}+4t^{19}+3t^{21}+2t^{23}+2t^{25}+t^{27}\\
1+t^2+t^4+2t^6+2t^8+3t^{10}+3t^{12}+2t^{14}\\\quad+3t^{16}+2t^{18}+2t^{20}+2t^{22}+t^{24}+t^{26}+t^{28}\\
t+t^5+t^7+t^9+2t^{11}+t^{13}\\\quad+t^{15}+2t^{17}+t^{19}+t^{21}+t^{23}+t^{27}\\
t^2+t^4+t^6+2t^8+2t^{10}+2t^{12}+2t^{14}\\\quad+2t^{16}+2t^{18}+2t^{20}+t^{22}+t^{24}+t^{26}
 \end{array}
\right)
\]

\[(H_A(t)_{i,7})_{1\leq i\leq8}=\left(
 \begin{array}{l}
 t^6+t^{12}+t^{16}+t^{22}\\
 t^5+t^7+t^{11}+t^{13}+t^{15}+t^{17}+t^{21}+t^{23}\\
 t^4+t^6+t^8+t^{10}+t^{12}+2t^{14}\\\quad+t^{16}+t^{18}+t^{20}+t^{22}+t^{24}\\
 t^3+t^5+t^7+2t^9+t^{11}+2t^{13}+2t^{15}\\\quad+t^{17}+2t^{19}+t^{21}+t^{23}+t^{25}\\
 t^2+t^4+t^6+2t^8+2t^{10}+2t^{12}+2t^{14}\\\quad+2t^{16}+2t^{18}+2t^{20}+t^{22}+t^{24}+t^{26}\\
 t+t^5+t^7+t^9+2t^{11}+t^{13}+t^{15}\\\quad+2t^{17}+t^{19}+t^{21}+t^{23}+t^{27}\\
 1+t^6+t^{10}+t^{12}+t^{16}+t^{18}+t^{22}+t^{28}\\
 t^3+t^7+t^9+t^{11}+t^{13}+t^{15}\\\quad+t^{17}+t^{19}+t^{21}+t^{25}
 \end{array}
\right)
\]

\[
(H_A(t)_{i,8})_{1\leq i\leq8}
=\left(
 \begin{array}{l}
t^5+t^9+t^{13}+t^{15}+t^{19}+t^{23}\\
t^4+t^6+t^8+t^{10}+t^{12}+2t^{14}\\\qquad+t^{16}+t^{18}+t^{20}+t^{22}+t^{24}\\
t^3+t^5+2t^7+t^9+2t^{11}+2t^{13}+2t^{15}\\\quad+2t^{17}+t^{19}+2t^{21}+t^{23}+t^{25}\\
t^2+t^4+2t^6+2t^8+2t^{10}+3t^{12}+2t^{14}\\\quad+3t^{16}+2t^{18}+2t^{20}+2t^{22}+t^{24}+t^{26}\\
t+t^3+2t^5+2t^7+3t^9+3t^{11}+3t^{13}+3t^{15}\\\quad+3t^{17}+3t^{19}+2t^{21}+2t^{23}+t^{25}+t^{27}\\
t^2+t^4+t^6+2t^8+2t^{10}+2t^{12}+2t^{14}\\\quad+2t^{16}+2t^{18}+2t^{20}+t^{22}+t^{24}+t^{26}\\
t^3+t^7+t^9+t^{11}+t^{13}+t^{15}\\\quad+t^{17}+t^{19}+t^{21}+t^{25}\\
1+t^4+t^6+t^8+2t^{10}+t^{12}+2t^{14}\\\quad+t^{16}+2t^{18}+t^{20}+t^{22}+t^{24}+t^{28}
 \end{array}
\right)
\] 

\section{$HH^0(A)=Z$}
From the Hilbert series \ref{center} we see that we have one (unique up to a constant factor) central element of degree $2m_i-2$ for each exponent $m_i<\frac{h}{2}$. We will denote a $\deg i(<h-2)$ central element by $z_i$. 

From \ref{topdegree} and from the Hilbert series we can also see that the top degree $(=\deg h-2)$ center is spanned by one element $\omega_i$ in each $e_iAe_i$, such that $\nu(i)=i$.

The $\omega_i\in L[h-2]$ are already given in section \ref{Frobenius}, and we will find the $z_i\in U[-2]$ for each Dynkin quiver separately.

\subsection{$Q=D_{n+1}$}

We define the nonzero elements
\begin{align*}
b_{i,0}&=e_i,\\
b_{i,j}&=a_i^*\ldots a_{i+j-1}^*a_{i+j-1}\ldots a_i\, (\text{where } 1\leq j\leq\min\{i-1,n-1-i\}),\\
c_{i,j}&=a_i^*\ldots a_{n-2}^*(a_{n-2}a_{n-2}^*)a_{n-2}\ldots a_i\, (\text{ where } 1\leq i\leq n-2,\, 1\leq j\leq i-1\\
c_{n-1,j}&=(a_{n-2}a_{n-2}^*)^j,\,1\leq j\leq n-2\\
c_i'&=a_i^*\ldots a_{n-2}^*a_{n-1}^*a_{n-1}(a_{n-2}a_{n-2}^*)^{i-1} a_{n-2}\ldots a_2,\,1\leq i\leq n-1\\
d_0&=e_n,\\
d_j&=(a_{n-1}a_n^*a_na_{n-1}^*)^j \text{ for } 1\leq j<\frac{n}{2},\\
d_0'&=e_{n+1}, \\
d_j'&=(a_na_{n-1}^*a_{n-1}a_n^*)^j \text{ for } 1\leq j<\frac{n}{2}
\end{align*}
and extend this notation for any other $j$, where $b_{i,j}$, $c_{i,j}$, $d_j$ and  $d_j'$ are zero.

The exponents $m_i$ are $1,3,\ldots ,2n-1,n$ and $h=2n$. From \ref{center} we get  the Hilbert series of $Z$, depending on the parity of $n$, since $r_+=n+1$ for $n$ odd and $r_+=n-1$ for $n$ even:

\begin{align*}
n\text{ odd: }&h_Z(t)=1+t^4+t^8+\ldots+t^{2n-6}+(n+1)t^{2n-2},\\
n\text{ even: }&h_Z(t)=1+t^4+t^8\ldots+t^{2n-4}+(n-1)t^{2n-2}.
\end{align*}

The central elements of degree $4j<2n-2$ are
\[
z_{4j}=\sum\limits_{i=2j+1}^{n-1-2j}b_{i,2j}+\sum\limits_{i=0}^{2j-1}c_{n-1-i,2j-i}+d_j+d_j'
\]
since the relations \\
$b_{2j+1,2j}a_{2j}=0,$\\
$a_ib_{i,2j}=b_{i+1,2j}a_i$ for $2j+1<i<n-1-2j$,\\
$a_{n-1-2j}b_{n-1-2j,2j}=c_{n-2j,1}a_{n-1-2j}$,\\
$a_{n-1-i}c_{n-1-i,2j-1-i}=c_{n-i,2j-i}a_{n-1-i}$ for $0<i<2j-1$\\
$a_{n-1}c_{n-1,2j}=d_ja_{n-1}$, $a_nc_{n-1,2j}=d_j'a_{n}$\\
hold.

The top degree central elements are $\omega_i=c_i'$ ($1\leq i\leq n-1$), and additionally $\omega_n=d_{\frac{n-1}{2}}$, $-\omega_{n+1}=d_{\frac{n-1}{2}}'$ if $n$ is odd.\\

We have the multiplication laws $b_{i,2j}b_{i,2k}=b_{i,2(j+k)}$, \\
$d_jd_k=d_{j+k}$,  $d_j'd_k'=d_{j+k}'$ and $c_{i,j}c_{i,k}=c_{i,j+k+n-i-1}$ which implies that \\
$c_{n-1-i,2j-1-i}c_{n-1-i,2k-1-i}=c_{n-1-i,2j-1+2k-1+n-(n-1-i)-1}=c_{n-1-i,2(j+k)-i}$, so for $j+k<\frac{n-1}{2}$ we get the following product:
\[
z_{4j}z_{4k}=z_{4(j+k)}.
\]
If $n$ is odd and $j+k=\frac{n-1}{2}$, the multiplication becomes

\[
z_{4j}z_{4k}=d_{\frac{n-1}{2}}+d_{\frac{n-1}{2}}'=\omega_{n}-\omega_{n+1}.
\]

\subsection{$Q=E_6$}
The Coxeter number is $h=12$, and the exponents $m_i<\frac{h}{2}=6$ are $1$, $4$, $5$, $r_+=2$. For the center, we get the following Hilbert series (from Corollary \ref{center}):
\[
h_Z(t)=1+t^6+t^8+2t^{10}.
\]

From the degrees, we see that the product of any two positive degree central elements is always $0$. The central elements are $z_0=1$, $z_6$, $z_8$, $\omega_3$ and $\omega_6$.

We use the notation $x_i=a_i^*a_i$ and give the central elements $z_6$ and $z_8$ explicitly:

\begin{proposition}
A central element of $\deg 6$ is given by
\[
z_6=a_1a_2x_3a_2^*a_1^*-a_2x_3^2a_2^*-x_5x_3x_5+a_3x_2^2a_3^*-a_4a_3x_2a_3^*a_4^*
\]
\end{proposition}
\begin{proof}
Observe that $z_6^*=z_6$, so it is enough to show that $z_6$ commutes with all $a_i$. Since $\eta(z_6)=z_6$, we have only to show that $z_6$ commutes with $a_1$, $a_2$, $a_5$:

\begin{align*}
a_1z_6&=-a_1a_2x_3^2a_2^*=-a_1a_2x_5x_2a_2^*=a_1a_2x_3a_2^*a_1^*a_1=z_6a_1,\\
a_2z_6&=-a_2x_5x_3x_5=-a_2x_3x_2x_5=-a_2x_3x_3x_2=z_6a_2,\\
a_5z_6&=a_5x_5x_3x_5=0=z_6a_5.
\end{align*}

\end{proof}

\begin{proposition}
The $\deg 8$ central element of $A$ is given by
\[
z_8=-a_2x_5x_3x_5a_2^*-x_5x_3^2x_5-a_3x_5x_2x_5a_3^*.
\]
\end{proposition}

\begin{proof}
Since $z_8^*=z_8$ and $\eta(z_8)=z_8$, it is enough to show that $z_8$ commutes with $a_1$, $a_2$ and $a_5$:

\begin{align*}
a_1z_8&=-a_1a_2x_5x_3x_5a_2^*=-a_1a_2x_5x_2x_3a_2^*=-a_1a_2x_2x_3x_3a_2^*=0=a_1z_8,\\
a_2z_8&=-a_2x_5x_3^2x_5=-a_2x_5x_2x_3x_2=-a_2x_5x_3x_5x_2=z_8a_2,\\
a_5z_8&=a_5x_5x_3x_5x_3=0=z_8a_5.
\end{align*}
\end{proof}

\subsection{$Q=E_7$}
The Coxeter number is $h=18$, the exponents $m_i<\frac{h}{2}=9$ are $1,5,7$, $r_+=7$, and the Hilbert series of the center is (see Corollary \ref{center}):
\[
h_Z(t)=1+t^8+t^{12}+7t^{16}
\]

The center is spanned by $z_0=1$, $z_8$, $z_{12}$, $\omega_1,\ldots,\omega_7$. The only interesting product to compute is $z_8^2$ which lies in the top degree.

We give $z_8$ and $z_{12}$ explicitly:
\begin{proposition}
The $\deg 8$ central element of $A$ is given by
\begin{align*}
z_8&=-a_1a_2a_3x_6a_3^*a_2^*a_1^*-a_2a_3x_4^2a_3^*a_2^*-a_3x_6x_4x_6a_3^*-x_4x_3^2x_4\\
&\quad-a_4x_4x_6x_4a_4^*+a_6x_4x_6x_4a_6^*.
\end{align*}

\end{proposition}
\begin{proof}
It is clear that $z_8^*=z_8$, so it is enough to prove that $z_8$ commutes with every $a_i$:
\begin{align*}
a_1z_8&=-a_1a_2a_3x_4^2a_3^*a_2^*=-a_1a_2a_3x_6x_3a_3^*a_2^*=-a_1a_2a_3x_6a_3^*a_2^*x_1=z_8a_1,\\
a_2z_8&=-a_2a_3x_6x_4x_6a_3^*=-a_2a_3x_4x_3x_6a_3^*=-a_2a_3x_4^2x_3a_3^*=-a_2a_3x_4^2a_3^*x_2= z_8a_2,\\
a_3z_8&=-a_3x_4x_3^2x_4=-a_3x_6x_4x_3x_4=-a_3x_6x_4x_6x_3=z_8a_3,\\
a_4z_8&=-a_4x_4x_3^2x_4=-a_4x_4x_6x_4^2=z_8a_4\\
a_5z_8&=-a_5a_4x_4x_6x_4a_4^*=0=z_8a_5 \\
a_6z_8&=-a_6x_4x_3^2x_4=-a_6x_3^3x_6=a_6x_4x_6x_4x_6=z_8a_6.
\end{align*}

\end{proof}

\begin{proposition}
The $\deg 12$ central element of $A$ is
\[
z_{12}=-a_3x_4x_6x_4x_6x_4a_3^*-x_4x_6x_4^2x_6x_4+a_4x_6x_4x_6x_4x_6a_4^*+a_6x_4x_6x_4x_6x_4a_6^*.
\]
\end{proposition}
\begin{proof}
$z_{12}^*=z_{12}$ is clear, so we have only to show that $z_{12}$ commutes with all $a_i$:

\begin{align*}
a_1z_{12}&=0=z_{12}a_1,\\
a_2z_{12}&=-a_2a_3x_4x_6x_4x_6x_4a_3^*=a_2a_3x_4x_6x_3^2x_6a_3^*=a_2a_3x_4x_6x_4^2x_6a_3^*\\
&=a_2a_3x_3^2x_4^2x_6a_3^*=0=z_{12}a_2,\\
a_3z_{12}&=-a_3x_4x_6x_4^2x_6x_4=a_3x_4x_6x_3^2x_6x_3=-a_3(x_4x_6)^2x_4x_3=z_{12}a_3,\\
a_4z_{12}&=-a_4x_4x_6x_4^2x_6x_4=a_4x_3^2x_4x_3x_6x_4=a_4(x_6x_4)^2x_6=z_{12}a_4,\\
a_5z_{12}&=a_5a_4x_6x_4x_6x_4x_6a_4^*=-a_5a_4x_6x_3^3x_6a_4^*=a_5a_4x_3^4x_6a_4^*=0=z_{12}a_5,\\
a_6z_{12}&=-a_6x_4x_6x_4^2x_6x_4=a_6x_4x_6x_3^2x_6x_3=-a_6x_4x_6x_4x_6x_4x_3=\\
&=a_6(x_4x_6)^3+a_6(x_4x_6)^2x_4^2=z_{12}a_6+a_6x_3^4x_4^2=z_{12}a_6.
\end{align*}
\end{proof}

Now we calculate the product $z_8^2$:
\begin{proposition}
\[z_8^2=\omega_1+\omega_3-\omega_7.\]
\end{proposition}
\begin{proof}
We get
\begin{align*}
z_8^2&=(a_1a_2a_3x_6a_3^*a_2^*a_1^*)^2+(a_2a_3x_4^2a_3^*a_2^*)^2+(a_3x_6x_4x_6a_3^*)^2+(x_4x_3^2x_4)^2\\
&\quad+(a_4x_4x_6x_4a_4^*)^2+(a_6x_4x_6x_4a_6^*)^2\\
&=\underbrace{a_1a_2a_3x_6x_3^3x_6a_3^*a_2^*a_1^*}_{=\omega_1}+\underbrace{a_2a_3x_4^2x_3^2x_4^2a_3^*a_2^*}_{=a_2a_3x_4^2x_6^2x_4^2a_3^*a_2^*=0}+\underbrace{a_3x_6x_4x_6x_3x_6x_4x_6a_3^*}_{=\omega_3}\\
&\quad+\underbrace{x_4x_3^2x_4^3x_3^2x_4}_{=0}+\underbrace{a_4x_4x_6x_4^3x_6x_4a_4^*}_{=0}+\underbrace{a_6(x_4x_6)^3x_4a_6^*}_{=-\omega_7}
\end{align*}

\end{proof}

\subsection{$Q=E_8$}
The Coxeter number $h=30$, and the exponents $m_i<\frac{h}{2}=15$ are $1,7,11,13$, $r_+=8$. For the center, we get the following Hilbert series (from Corollary \ref{center}):
\[
h_Z(t)=1+t^{12}+t^{20}+t^{24}+8t^{28}.
\]

The center is spanned by $z_0=1$, $z_{12}$, $z_{20}$, $z_{24}$, $\omega_1,\ldots,\omega_8$. The only interesting product is $z_{12}^2$.

\begin{proposition}
The $\deg 12$ central element of $A$ is
\begin{align*}
z_{12}&=a_1a_2a_3x_6x_4x_6a_3^*a_2^*a_1^*+a_2a_3x_4x_3^2x_4a_3^*a_2^*+a_3(x_4x_6)^2x_4a_3^*\\
&\quad+(x_3x_4x_3)^2-a_4(x_6x_4)^2x_6a_4^*+a_5a_4x_6x_4^2x_6a_4^*a_5^*-a_6(x_4x_6)^2x_4a_6^*.
\end{align*}
\end{proposition}
\begin{proof}
$z_{12}^*=z_{12}$, so we only have to show that $z_{12}$ commutes with all $a_i$:
\begin{align*}
a_0z_{12}&=a_0a_1a_2a_3x_6x_4x_6a_3^*a_2^*a_1^*=a_0a_1a_2a_3x_6x_3x_4a_3^*a_2^*a_1^*\\&=a_0a_1a_2a_3x_3x_4^2a_3^*a_2^*a_1^*=0=z_{12}a_0,\\
a_1z_{12}&=a_1a_2a_3x_4x_3^2x_4a_3^*a_2^*=a_1a_2a_3x_6x_4x_6x_3a_3^*a_2^*\\&=a_1a_2a_3x_6x_4x_6a_3^*a_2^*x_1=z_{12}a_1,\\
a_2z_{12}&=a_2a_3(x_4x_6)^2x_4a_3^*=-a_2a_3x_4x_3^2x_6x_3a_3^*=a_2a_3x_4x_3^2x_4a_3^*x_2=z_{12}a_2,\\
a_3z_{12}&=a_3(x_3x_4x_3)^2=-a_3x_3x_6x_3^2x_4x_3=a_3x_4x_6x_4x_6x_4x_3=z_{12}a_3,\\
a_4z_{12}&=a_4(x_3x_4x_3)^2=-a_4(x_6x_4)^3=z_{12}a_4,\\
a_5z_{12}&=-a_5a_4(x_6x_4)^2x_6a_4^*=a_5a_4x_3^4x_6a_4^*=-a_5a_4x_3^4x_4a_4^*\\
&=a_5a_4x_6x_4^2x_6a_4^*x_5=z_{12}a_5,\\
a_6z_{12}&=a_6(x_3x_4x_3)^2=-a_6(x_6x_4)^2x_6=z_{12}a_6.
\end{align*}
\end{proof}

\begin{proposition}
The central element of $\deg 20$ is
\begin{align*}
z_{20}&=-a_1a_2a_3x_4^2x_3^3x_4^2a_3^*a_2^*a_1^*-a_2a_3(x_6x_4)^2(x_4x_6)^2a_3^*a_2^*+a_3(x_6x_4)^4x_6a_3^*\\
&\quad-(x_4x_6)^5+(x_6x_4^2)^3x_6-(x_6x_4)^5\\
&\quad-a_4(x_4x_6x_4)^3a_4^*-a_6(x_4x_6)^4x_4a_6^*.
\end{align*}
\end{proposition}

\begin{proof}
$z_{20}^*=z_{20}$, so it is sufficient to show that $z_{20}$ commutes with all $a_i$:
\begin{align*}
a_0z_{20}&=-a_0a_1a_2a_3x_4^2x_3^3x_4^2a_3^*a_2^*a_1^*=a_0a_1a_2a_3x_4^2x_6x_4x_6x_4^2a_3^*a_2^*a_1^*\\
&=-a_0a_1a_2a_3x_4^2x_6x_4x_6x_4x_6a_3^*a_2^*a_1^*=-a_0a_1a_2a_3x_4^2x_3^4x_6a_3^*a_2^*a_1^*\\
&=-a_0a_1a_2a_3x_6x_3^5x_6a_3^*a_2^*a_1^*=0=z_{20}a_0,\\
a_1z_{20}&=-a_1a_2a_3(x_6x_4)^2(x_4x_6)^2a_3^*a_2^*=-a_1a_2a_3x_4x_3x_6x_4(x_4x_6)^2a_3^*a_2^*\\
&=-a_1a_2a_3x_4^2x_3^3x_4x_3x_6a_3^*a_2^*=-a_1a_2a_3x_4^2x_3^3x_4^2x_3a_3^*a_2^*\\
&=-a_1a_2a_3x_4^2x_3^3x_4^2a_3^*a_2^*x_1=z_{20}a_1,\\
a_2z_{20}&=a_2a_3(x_6x_4)^4x_6a_3^*=a_2a_3x_6x_4x_6x_3^3x_6x_3x_6a_3^*\\
&=-a_2a_3x_6x_4x_6x_3^3x_4^2x_3a_3^*=a_2a_3(x_6x_4)^2x_4x_6x_4^2x_3a_3^*\\
&=-a_2a_3(x_6x_4)^2(x_4x_6)^2a_3^*x_2=z_{20}a_2.\\
\end{align*}
From
\[
 a_3(x_6x_4^2)^3x_6=a_3x_4x_3^3(x_4^2x_6)^2=a_3x_4x_6x_4x_6x_3^2x_4x_3^2x_6=a_3(x_4x_6)^5
\]
follows
\[
 a_3z_{20}=a_3(-(x_4x_6)^5+(x_6x_4^2)^3x_6-(x_6x_4)^5)=a_3(x_6x_4)^4x_6x_3=z_{20}a_3.
\]
To show that $z_{20}$ commutes with $a_4$, we use the following rules:
\begin{align*}
a_4(x_4x_6)^3&=a_4x_4x_3^4x_6\\
&=-a_4x_6x_4^2x_6x_4x_6-a_4x_6x_4x_6x_4^2x_6\\
&=-a_4(x_4x_6x_4)^2-a_4(x_4x_6)^2x_4^2
&=a_4(x_6x_4^2)^2+a_4(x_6x_4)^2,
\end{align*}
so
\begin{align*}
a_4(x_4x_6)^5&=-a_4x_6x_4^2x_6(x_4x_6)^3-a_4(x_6x_4x_6x_4)^2(x_4x_6)^3\\
&=a_4x_6x_4 (x_4x_6x_4)^2 x_4x_6+a_4x_6x_4(x_4x_6)^2x_4^3x_6\\
&\quad-a_4(x_6x_4)^2(x_6x_4^2)^2-a_4(x_6x_4)^5,
\end{align*}
\begin{align*}
a_4z_{20}&=-a_4(x_4x_6)^5+a_4(x_6x_4^2)^3x_6-a_4(x_6x_4)^5=a_4(x_6x_4)^2(x_6x_4^2)^2\\
&=a_4x_6x_3^4x_4^2x_6x_4^2=a_4(x_4x_6x_4)^3x_4=z_{20}a_4,\\
a_5z_{20}&=-a_5a_4(x_4x_6x_4)^3a_4^*=0=z_{20}a_5,\\
a_6z_{20}&=-a_6(x_4x_6)^5+a_6(x_6x_4^2)^3x_6-a_6(x_6x_4)^5=-a_6(x_4x_6)^5=z_{20}a_6.
\end{align*}
\end{proof}

In $\deg 24$ there is exactly one central element (up to constant multiple). Since we already have a $\deg 12$ central element, the most obvious candidate is its square. In fact, this is true from the following

\begin{proposition} 
The $\deg 24$ central element of $A$ is 
\[
z_{24}=z_{12}^2.
\]
\end{proposition}
\begin{proof}
A simple calculation shows that $z_{12}^2$ is nonzero:
\[
x_4x_3z_{12}^2=\omega_4.
\]
\end{proof}

\section{$HH^1(A)$}
Recall Theorem \ref{t1} where we know that $HH^1(A)$ is isomorphic to the non-topdegree part of $HH^0(A)$. In fact, $HH^1(A)$ is generated by the central elements in the following way:
\begin{proposition}
$HH^1(A)$ is spanned by maps
\[
\begin{array}{c}
\theta_k:(A\otimes V\otimes A)\longrightarrow A,\\
\theta_k (1\otimes a_i \otimes 1)=0,\\
\theta_k (1\otimes a_i^*\otimes 1)=a_i^*z_k.
\end{array}
\]
\end{proposition}
\begin{proof}
These maps clearly lie in $\ker d_2^*$: Recall
\[
\begin{array}{ccc}
A\otimes A&\stackrel{d_2}{\longrightarrow}& A\otimes V\otimes A\\
x\otimes y&\longmapsto & \sum\limits_{a\in\bar Q}\epsilon_axa\otimes a^*\otimes y+\sum\limits_{a\in\bar Q}\epsilon_ax\otimes a\otimes a^*y,
\end{array}
\]
then
\begin{align*}
 d_2^*\circ\theta_k(1\otimes1)&=\theta_k(\sum\limits_{a\in\bar Q}\epsilon_aa\otimes a^*\otimes 1+\sum\limits_{a\in\bar Q}\epsilon_a1\otimes a\otimes a^*)\\
 &=\sum\limits_{i\in I}a_ia_i^*z_k-\sum\limits_{i\in I}a_i^*a_iz_k=\sum\limits_{i\in I}[a_i,a_i^*]z_k=0.
\end{align*}

We will later see in section \ref{HH4} that $HH^4(A)$ is generated by $\zeta_k$ where $\zeta_k(\theta_k)=1$ under the duality $HH^4(A)=(HH^1(A))^*$ established in \cite{EE2}, so $\theta_k$ is nonzero in $HH^1(A)$.
\end{proof}

\section{$HH^2(A)$}

We know from \ref{t1} that $HH^2(A)=K[-2]$ lies in degree $-2$, i.e. in the lowest degree of $A^R[-2]$ (using the identifications in \cite[Section 4.5.]{EE2}), that is in $R[-2]$. Since the image of $d_2^*$ lies in degree $>-2$, $HH^2(A)=\ker d_3^*$.

\begin{proposition}
$HH^2(A)$ is given by the kernel of the matrix $H_A(1)$, where we identify $\mathbb{C}^I=R=\oplus_{i\in I}Re_i$.
\end{proposition}
\begin{proof}
Recall

\[
 d_3^*(x)=\sum\limits_{x_i\in B} x_i x x_i^*=\sum\limits_{j,k\in I}\sum\limits_{x_i\in B_{j,k}} x_ixx_i^*.
\]
For each $x_i\in e_kAe_j$, we see that $x_ie_lx_i^*=\delta_{jl}\omega_k$.

It follows that for $x=\sum\limits_{i\in I} \lambda_ie_i$ the map is given by

\[
 d_3^*(x)=\sum\limits_{i\in I}\mu_i\omega_i,
\]
where the vectors $\lambda=(\lambda_i)_{i\in I}\in \mathbb{C}^I$ and $\mu=(\mu_i)_{i\in I}\in\mathbb{C}^I$ satisfy the equation

\begin{equation}\label{d3-star}
 H_A(1)\lambda=\mu.
\end{equation}
So the kernel of $d_3^*$ is given by the kernel of $H_A(1)$.
\end{proof}
Now, we find the elements in $HH^2(A)$ for the quivers separately.
\subsection{$Q=D_{n+1}$, $n$ even}
\begin{equation}\label{Hilbert D}
 H_A(1)=\left(
 \begin{array}{ccccccccc}
  2&2&2&2&\ldots&\ldots&2&1&1\\
  2&4&4&4&\ldots&\ldots&4&2&2\\
  2&4&6&6&\ldots&\ldots&6&3&3\\
  2&4&6&8&\ldots&\ldots&8&4&4\\
  \vdots&\vdots&\vdots&\vdots&\ddots&&\vdots&\vdots&\vdots\\
  \vdots&\vdots&\vdots&\vdots&&\ddots&\vdots&\vdots&\vdots\\
  2&4&6&8&\ldots&\ldots&2(n-1)&n-1&n-1\\
  1&2&3&4&\ldots&\ldots&n-1&\frac{n}{2}&\frac{n}{2}\\
  1&2&3&4&\ldots&\ldots&n-1&\frac{n}{2}&\frac{n}{2}\\  
 \end{array}
 \right)
\end{equation}
with kernel $\langle e_n-e_{n+1}\rangle$.
So a basis of $HH^2(A)$ is given by 
\[\{f_n=[e_n-e_{n+1}]\}.\]

\subsection{$Q=E_6$}

\begin{equation}\label{Hilbert E}
 \left(
 \begin{array}{cccccc}
  2&3&4&3&2&2\\
  3&6&8&6&3&4\\
  4&8&12&8&4&6\\
  3&6&8&6&3&4\\
  2&3&4&3&2&2\\
  2&4&6&4&2&4
 \end{array}
 \right)
\end{equation}
with kernel $\langle e_1-e_5, e_2-e_4\rangle$. So a basis of $HH^2(A)$ is given by \[\{f_1=[e_1-e_5],f_2=[e_2-e_4]\}.\]
\section{$HH^3(A)$}
We know that $HH^3(A)$ lives in degree $-2$. Using the notations and identifications in \cite[Section 4.5.]{EE2}, we see that the kernel of $d_4^*$ has to be the top degree part of $\mathfrak{N}^R[-h]$ (since $Im\,d_3^*$ lives in degree $-2$), so
\[
 HH^3(A)=\mathfrak{N}^R[-h](-2)/Im\,d_3^*.
\]
\begin{proposition}
  $HH^3(A)$ is given by the cokernel of the matrix $H_A(1)$, where we identify $\mathbb{C}^I=A^{top}=\oplus_{i\in I}e_iA^{top}e_{\nu(i)}$.
\end{proposition}
\begin{proof}
This follows immediately from the discussion in the previous section because $d_3^*$ is given by $H_A(1)$.
\end{proof}

Note that $HH^3(A)=(HH^2(A))^*$ under the duality in \cite{EE2}. We choose a basis $h_i$ of $HH^3(A)$, so that $h_i(f_j)=\delta_{i,j}$
\subsection{$Q=D_{n+1}$, $n$ even}
From $H_A(1)$ in \ref{Hilbert D} we see that:
\begin{eqnarray*}
d_3^*(2e_1-e_2)&=&2\omega_1,\\
d_3^*(-e_{i-1}+2e_i-e_{i+1})&=&2\omega_i\qquad\forall 2\leq i\leq n-2,\\
d_3^*((-n-1)e_{n-2}+2(n-1)e_{n-1}-2(n-1)e_n)&=&(n-1)\omega_{n-1},\\
d_3^*(2e_n-e_{n-1})&=&\omega_n+\omega_{n+1},
\end{eqnarray*}
so \[HH^3(A)=(\mathfrak{N}^R)^{top}[-h]/(\omega_1=\omega_2=\ldots=\omega_{n-1}=0,\,\omega_n+\omega_{n+1}=0)\]
with basis \[\{h_n=[\omega_n]\}\]

\subsection{$Q=E_6$}
From $H_A(1)$ in \ref{Hilbert E} we see that:

\begin{eqnarray*}
d_3^*(2e_1-e_2)&=&\omega_1+\omega_5,\\
d_3^*(-e_1+2e_2-e_3)&=&\omega_2+\omega_4,\\
d_3^*(-2e_2+2e_3-e_6)&=&2\omega_3,\\
d_3^*(-e_3+2e_6)&=&2\omega_6,
\end{eqnarray*}
so 
\[HH^3(A)=(\mathfrak{N}^R)^{top}[-h]/(\omega_3=\omega_6=\omega_1+\omega_5=\omega_2+\omega_3=0)\]
with basis \[\{h_1=[\omega_1],\,h_2=[\omega_2]\}.\]
\section{$HH^4(A)$}\label{HH4}
We have $HH^4(A)=U^*[-2]$, so its top degree is $-4$, and its generators sit in degrees $-4-\deg z_k$ for each central element, one in each degree.
\begin{proposition}
 Let $\zeta_0\in\ker d_5^*$ be a top degree element in\\ $(V\otimes\mathfrak{N})^R[-h-2]$, such that $m(\zeta_0)$ is nonzero, where $m$ is the multiplication map. Then $HH^4(A)$ is generated by elements $\zeta_k\in\ker d_5^*$ which satisfy $\zeta_kz_k=\zeta_0$.
\end{proposition}
\begin{proof}
 If $x\in \mathfrak{N}^R[-h])$ lies in degree $-4$, then $m(d_4^*(x))=0$, so  $\zeta_0$ is nonzero in $HH^4(A)$. 

For every non-topdegree central element $z_k$ we can find a $\zeta_k$ satisfying the properties above, which is done for each quiver separately below.

For any central element $z\in A$, we have that $d_4^*(zy)=d_4^*(y)z$. If  $\zeta_k=d_4^*(y)$, then by construction $\zeta_0=\zeta_kz_k=d_4^*(z_ky)$ which is a contradiction.
 
So these $\zeta_k$ are all nonzero in $HH^4(A)$, and also generate this cohomology space.

\end{proof}

A basis of $HH^4(A)$ is given by these $\zeta_k$, and we choose them so that $\zeta_k(\theta_k)=1$ under the duality $HH^4(A)=(HH^1(A))^*$ in \cite{EE2}.

\subsection{$Q=D_{n+1}$, $n$ odd}
We define
\begin{align*}
 \zeta_0=&[a_{n-1}^*\otimes a_{n-1}a_n^*a_n(a_{n-1}^*a_{n-1}a_n^*a_n)^{\frac{n-3}{2}}+
 &a_{n-1}\otimes a_{n}^*a_na_{n-1}^*(a_{n-1}a_n^*a_na_{n-1}^*)^{\frac{n-3}{2}}],\\
 \zeta_{4k}=
 &\frac{1}{2}[a_{n-1}^*\otimes a_{n-1}a_n^*a_n(a_{n-1}^*a_{n-1}a_n^*a_n)^{\frac{n-3}{2}-k}+
 &a_{n-1}\otimes a_{n}^*a_na_{n-1}^*(a_{n-1}a_n^*a_na_{n-1}^*)^{\frac{n-3}{2}-k}-\\
 &a_{n}^*\otimes a_{n}a_{n-1}^*a_{n-1}(a_{n}^*a_{n}a_{n-1}^*a_{n-1})^{\frac{n-3}{2}-k}-
 &a_{n}\otimes a_{n-1}^*a_{n-1}a_{n}^*(a_{n}a_{n-1}^*a_{n-1}a_{n}^*)^{\frac{n-3}{2}-k}].
\end{align*}
\subsection{$Q=D_{n+1}$, $n$ even}
We define
\begin{align*}
 \zeta_0=&[a_{n-1}^*\otimes a_{n-1}(a_n^*a_na_{n-1}^*a_{n-1})^{\frac{n-2}{2}-k}
 &+a_{n-1}\otimes a_{n}^*(a_na_{n-1}^*a_{n-1}a_n^*)^{\frac{n-2}{2}-k}],\\
 \zeta_{4k}=
 &\frac{1}{2}[a_{n-1}^*\otimes a_{n-1}(a_n^*a_na_{n-1}^*a_{n-1})^{\frac{n-2}{2}-k}
 &+a_{n-1}\otimes a_{n}^*(a_na_{n-1}^*a_{n-1}a_n^*)^{\frac{n-2}{2}-k}-\\
 &a_{n}^*\otimes a_{n}(a_{n-1}^*a_{n-1}a_{n}^*a_{n})^{\frac{n-2}{2}-k}-
 &a_{n}\otimes a_{n-1}^*(a_{n-1}a_{n}^*a_{n}a_{n-1}^*)^{\frac{n-2}{2}-k}].
\end{align*}
\subsection{$Q=E_6$}
We define
\begin{eqnarray*}
 \zeta_0&=&[a_3^*\otimes a_3(a_2^*a_2a_3^*a_3)^2+a_3\otimes a_2^*(a_2a_3^*a_3a_2^*)^2],\\
 \zeta_6&=&\frac{1}{4}[-a_3^*\otimes a_3a_2^*a_2-a_3\otimes a_2^*a_2a_2^*+a_2^*\otimes a_2a_2^*a_2+a_2\otimes a_2^*a_2a_3^*\\
 &&\qquad -a_2^*\otimes a_2a_3^*a_3-a_2\otimes a_3^*a_3a_3^*+a_3^*\otimes a_3a_3^*a_3+a_3\otimes a_3^*a_3a_2^*],\\
 \zeta_8&=&\frac{1}{2}[a_3^*\otimes a_3+a_3\otimes a_2^*-a_2^*\otimes a_2-a_2\otimes a_3^*].
\end{eqnarray*}
\subsection{$Q=E_7$}
We define
\begin{eqnarray*}
 \zeta_0&=&[a_4^*\otimes a_4a_3^*a_3(a_4^*a_4a_3^*a_3)^3+a_4\otimes a_3^*a_3a_4^*(a_4a_3^*a_3a_4^*)^3],\\
 \zeta_8&=&\frac{1}{2}[a_4^*\otimes a_4a_3^*a_3a_4^*a_4a_3^*a_3+a_4\otimes a_3^*a_3a_4^*a_4a_3^*a_3a_4^*\\
 &&
 -a_3^*\otimes a_3a_4^*a_4a_3^*a_3a_4^*a_4-a_3\otimes a_4^*a_4a_3^*a_3a_4^*a_4a_3^*],\\
 \zeta_{12}&=&\frac{1}{2}[a_4^*\otimes a_4a_3^*a_3+a_4\otimes a_3^*a_3a_4^*-a_3^*\otimes a_3a_4^*a_4-a_3\otimes a_4^*a_4a_3^*].
\end{eqnarray*}

\subsection{$Q=E_8$}
We define
\begin{eqnarray*}
 \zeta_0&=&[a_4^*\otimes a_4a_3^*a_3(a_4^*a_4a_3^*a_3)^6+a_4\otimes a_3^*a_3a_4^*(a_4a_3^*a_3a_4^*)^6],\\
 \zeta_{12}&=&\frac{1}{2}[a_4^*\otimes a_4a_3^*a_3(a_4^*a_4a_3^*a_3)^3+a_4\otimes a_3^*a_3a_4^*(a_4a_3^*a_3a_4^*)^3\\
 &&-a_3^*\otimes a_3a_4^*a_4(a_3^*a_3a_4^*a_4)^3-a_3\otimes a_4^*a_4a_3^*(a_3a_4^*a_4a_3^*)^3],\\
 \zeta_{20}&=&\frac{1}{2}[a_4^*\otimes a_4a_3^*a_3a_4^*a_4a_3^*a_3+a_4\otimes a_3^*a_3a_4^*a_4a_3^*a_3a_4^*\\
 &&-a_3^*\otimes a_3a_4^*a_4a_3^*a_3a_4^*a_4-a_3\otimes a_4^*a_4a_3^*a_3a_4^*a_4a_3^*],\\
 \zeta_{24}&=&\frac{1}{2}[a_4^*\otimes a_4a_3^*a_3+a_4\otimes a_3^*a_3a_4^*-a_3^*\otimes a_3a_4^*a_4-a_3\otimes a_4^*a_4a_3^*].
\end{eqnarray*}

\section{$HH^5(A)$}

We have $HH^5(A)=U^*[-2]\oplus Y^*[-h-2]$. We discuss these two subspaces separately.
\subsection{$U^*[-2]$}\label{Ustar}
In $U^*[-2]$, like in $HH^4(A)$, we have generators coming from the center in some dual sense. 

We have $d_6^*(U^*[-2])=0$.

\begin{proposition}
 Let $\psi_0$ be a top degree element $[\omega_i]$ in some\\ $e_i\mathfrak{N}^Re_i[-h-2]$. Then $HH^5(A)$ is generated by $\psi_k\in\mathfrak{N}^R$ which satisfy $\psi_kz_k=\psi_0$.
\end{proposition}
\begin{proof}

If $\sum\limits_{a\in\bar Q} a\otimes x_a\in V\otimes\mathfrak{N}^R$ lies in degree $-4$, then the image of $d_5^*(x)=\sum\limits_{a} ax_a-x_a\eta(a)$, under the linear map $f$ (which is associated to $A$ as a Frobenius algebra) is zero where $f(\omega_i)=1$. So $\psi_0$ is nonzero in $HH^5(A)$.

For every non-topdegree central element $z_k$ we can find a $\zeta_k$ satisfying the properties above, which is done for each quiver separately in subsection \ref{HH5 results}.

For any central element $z\in A$, we have that $d_5^*(zy)=d_5^*(y)z$. If  $\psi_k=d_5^*(y)$, then by construction $\psi_0=\psi_kz_k=d_4^*(z_ky)$ which is a contradiction.

So these $\psi_k$ are nonzero in $HH^5(A)$ and generate this cohomology space.
\end{proof}
The relation $ax_a=x_a\eta(a)$ then gives us that all $\omega_i$'s are equivalent in $HH^5(A)$.
\subsection{$Y^*[-h-2]$}\label{Ystar}
We have to introduce some new notations.
\begin{definition}
We define $F$ to be the set of vertices in $I$ which are fixed by $\nu(i)$, i.e. \[F=\{i\in I|\nu(i)=i\}.\]
\end{definition}
\begin{definition}
 Let $\eta_{ij}$ be the restriction of $\eta$ on $e_iAe_j$ ($i,j\in F$). Let $n_{ij}^+=\dim\ker(\eta_{ij}-1)$ and $n_{ij}^-=\dim\ker(\eta_{ij}+1)$. 

We define the \emph{signed truncated dimension matrix} $(H_A^\eta)_{i,j\in F}$ in the following way:
\[
 (H_A^\eta)_{ij}=n_{ij}^+-n_{ij}^-.
\]
\end{definition}
Now we can make the following statement:
\begin{proposition}
 $Y^*[-h-2]$ is given by the kernel of the matrix $H_A^\eta$, where we identify $\mathbb{C}^F=\oplus_{i\in F}Re_i$.
\end{proposition}
\begin{proof}
$Y^*[-h-2]$ is the kernel of the restriction $d_6^*|_{\mathfrak{N}^R[-h-2](-h-2)=R_F[-h-2]}\rightarrow A^R[-2h]$, where $R_F$ is the linear span of $e_i$'s, such that $i$ is fixed by $\nu$,
\[
 d_6^*(x)=\sum\limits_{x_j\in B} x_jx\eta(x_j^*)=\sum\limits_{x_j\in B} \eta(x_j)xx_j^*.
\]

then \[d_6^*:R_F[-h-2]\rightarrow (A^{top})^R[-2h]\] can also be written as a matrix multiplication
\[
 H_A^\eta:\mathbb{C}^F\rightarrow\mathbb{C}^F
\]
under the identifications $R_F=\mathbb{C}^F=\bigoplus\limits_{i\in F}e_iA^{top}e_i$. 
\end{proof}
We compute the matrices $H_A^\eta$ and their kernels for each quiver separately.

Recall that $\dim Y=r_+-r_--\#\{m_i|m_i=\frac{h}{2}\}=\dim R_F-\#\{m_i|m_i=\frac{h}{2}\}$. We will find $Y^*$ explicitely for each quiver.

\subsubsection{$Q=E_6,E_8$}
$\frac{h}{2}$ is not an exponent, so $Y^*=R_F$.
\subsubsection{$Q=D_{n+1}$, $n$ odd}
All basis elements of $e_kAe_j$ given in section \ref{basis D_{n+1}} are eigenvectors of $\eta_{kj}$.

For any of these basis elements $x$, $\eta(x)=(-1)^{n_x}x$ where $n_x$ is the number of no-star letters in the monomial expression of $x$.
So $H_A^\eta$ can be computed directly.

For $k\leq j\leq n-1:$
\begin{eqnarray*}
 (H_A^\eta)_{k,j}&=&\left\{
 \begin{array}{lc}
\sum\limits_{l=0}^{k-1}(-1)^l+(-1)^{n-j}\sum\limits_{l=0}^{k-1}(-1)^l & k\leq n-j\\
\begin{array}{l}\sum\limits_{l=0}^{n-j-1}(-1)^l+(-1)^{n-j}\sum\limits_{l=0}^{k-1}(-1)^l\\
\quad+(-1)^{n-j}\sum\limits_{l=0}^{k-1+j-n}(-1)^l\end{array} & k>n-j
\end{array}
\right.\\
&=&\left\{
\begin{array}{lc}
 2 & k,l\,odd\\
 0 & else 
\end{array}
\right..
 \end{eqnarray*}
 
For $j<k\leq n-1:$
\begin{eqnarray*}
 (H_A^\eta)_{k,j}&=&\left\{
 \begin{array}{lc}
(-1)^{k-j}\sum\limits_{l=0}^{j-1}(-1)^l+(-1)^{n-j}\sum\limits_{l=0}^{j-1}(-1)^l & j\leq n-k\\
\begin{array}{l}(-1)^{k-j}\sum\limits_{l=0}^{n-k-1}(-1)^l+(-1)^{n-j}\sum\limits_{l=0}^{j-1}(-1)^l \\
\quad+(-1)^{n-j}\sum\limits_{l=0}^{k-1+j-n}(-1)^l\end{array} & j>n-k
\end{array}
\right.\\
&=&\left\{
\begin{array}{lc}
 2 & k,l\,odd\\
 0 & else 
\end{array}
\right..
 \end{eqnarray*}

For $k,j\leq n-1$:
\begin{eqnarray*}
(H_A^\eta)_{k,n}&=&\sum\limits_{l=0}^{k-1} (-1)^l=\left\{
\begin{array}{lc}
 1 & k\,odd\\
 0 & k\,even
\end{array}
\right.,\\
(H_A^\eta)_{n,n}&=&\sum\limits_{l=0}^{\frac{n-1}{2}}(-1)^{2l}=\frac{n+1}{2},\\
(H_A^\eta)_{n+1,n}&=&\sum\limits_{l=0}^{\frac{n-3}{2}}(-1)^{1+2l}=-\frac{n-1}{2},\\
(H_A^\eta)_{n,j}&=&(-1)^{n-j}\sum\limits_{l=0}^{j-1}(-1)^{l}=\left\{
\begin{array}{lc}
 1 & l\,odd\\
 0 & l\,even
\end{array}
\right..
\end{eqnarray*}

It is clear that $(H_A^\eta)_{k,n+1}=(H_A^\eta)_{k,n}$, $(H_A^\eta)_{n+1,n+1}=(H_A^\eta)_{n,n}$, \\$(H_A^\eta)_{n,n+1}=(H_A^\eta)_{n+1,n}$, $(H_A^\eta)_{n+1,j}=(H_A^\eta)_{n,j}$.

So the matrix $H_A^\eta$ is
\[
 H_A^\eta=
 \left(
 \begin{array}{ccccccc}
  2&0&\cdots&2&0&1&1\\
  0&0&\cdots&0&0&0&0\\
  \vdots&\vdots&\ddots&\vdots&\vdots&\vdots&\vdots\\
  2&0&\cdots&2&0&1&1\\
  0&0&\cdots&0&0&0&0\\
  1&0&\cdots&1&0&\frac{n+1}{2}&-\frac{n-1}{2}\\
  1&0&\cdots&1&0&-\frac{n+1}{2}&\frac{n+1}{2}
 \end{array}
 \right),
\]
and the kernel is given by
\[
 \langle e_{2k-1}-e_1,\, e_{2k}, (e_n+e_{n+1})-e_1|k\leq \frac{n-1}{2}\rangle.
\]

\subsubsection{$Q=D_{n+1}$, $n$ even}\label{D_{n+1}, n even}
Since $F=\{1,\ldots, n-1\}$, we work only with $e_kAe_j$ for $j,k\leq n-1$, and we have to work with a modified basis, so that they are all eigenvectors of $\eta$:

For \underline{$k\leq j\leq n-1$},
\begin{eqnarray*}
B_{k,j}&=&\{(a_{k-1}a_{k-1}^*)^la_k^*\cdots a_{j-1}^*|0\leq l\leq\min\{k-1,n-j-1\}\}\cup\\
&&\{(a_{k-1}a_{k-1}^*)^la_k^*\cdots (a_{n-1}^*a_{n-1}-a_n^*a_n)a_{n-2}a_j|0\leq l\leq k-1\}\cup\\
&&\{(a_{k-1}a_{k-1}^*)^la_k^*\cdots (a_{n-1}^*a_{n-1}+a_{n}^*a_{n})a_{n-2}a_j|0\leq l\leq k-1+j-n\}.
\end{eqnarray*}
For \underline{$j<k\leq n-1$},
\begin{eqnarray*}
B_{k,j}&=&\{a_{k-1}\cdots a_j(a_j^*a_j)^l|0\leq l\leq\min\{n-k-1,j-1\}\}\cup\\
&&\{a_k^*\cdots a_{n-2}^*(a_{n-1}^*a_{n-1}-a_n^*a_n)a_{n-2}\cdots a_j(a_j^*a_j)^l|0\leq l\leq j-1\}\cup\\
&&\{a_k^*\cdots a_{n-2}^*(a_{n-1}^*a_{n-1}+a_{n}^*a_{n})a_{n-2}\cdots a_j(a_j^*a_j)^l|0\leq l\leq j-1+k-n\}.
\end{eqnarray*}

For $k\leq j\leq n-1:$
\begin{eqnarray*}
 (H_A^\eta)_{k,j}&=&\left\{
 \begin{array}{lc}
\sum\limits_{l=0}^{k-1}(-1)^l+(-1)^{n-j-1}\sum\limits_{l=0}^{k-1}(-1)^l & k\leq n-j\\
\begin{array}{l}\sum\limits_{l=0}^{n-j-1}(-1)^l+(-1)^{n-j-1}\sum\limits_{l=0}^{k-1}(-1)^l\\
\quad+(-1)^{n-j}\sum\limits_{l=0}^{k-1+j-n}(-1)^l\end{array} & k>n-j
\end{array}
\right.\\
&=&\left\{
\begin{array}{lc}
 2 & k,l\,odd\\
 0 & else 
\end{array}
\right..
 \end{eqnarray*}
 
For $j<k\leq n-1:$
\begin{eqnarray*}
 (H_A^\eta)_{k,j}&=&\left\{
 \begin{array}{lc}
(-1)^{k-j}\sum\limits_{l=0}^{j-1}(-1)^l+(-1)^{n-j-1}\sum\limits_{l=0}^{j-1}(-1)^l & j\leq n-k\\
\begin{array}{l}(-1)^{k-j}\sum\limits_{l=0}^{n-k-1}(-1)^l+(-1)^{n-j-1}\sum\limits_{l=0}^{j-1}(-1)^l \\
\quad+(-1)^{n-j}\sum\limits_{l=0}^{k-1+j-n}(-1)^l\end{array} & j>n-k
\end{array}
\right.\\
&=&\left\{
\begin{array}{lc}
 2 & k,l\,odd\\
 0 & else 
\end{array}
\right..
 \end{eqnarray*}

So the matrix $H_A^\eta$ is
\[
 H_A^\eta=
 \left(
 \begin{array}{cccccc}
 2&0&\cdots&2&0&2\\
 0&0&\cdots&0&0&0\\
 \vdots&\vdots&\ddots&\vdots&\vdots&\vdots\\
 2&0&\cdots&2&0&2\\
 0&0&\cdots&0&0&0\\
 2&0&\cdots&2&0&2
 \end{array}
 \right),
\]
and we get immediately its kernel
\[
 \langle e_{2k+1}-e_1,e_{2k}|1\leq k\leq \frac{n}{2}\rangle.
\]

\subsubsection{$Q=E_7$}
We don't use an explicit basis of $A$ here. All we have to know is the number of no-star letters in the monomial basis elements which can be directly obtained from the Hilbert series $H_A(t)$ in the following way:
given a monomial $x$ of length $l$ in $e_kAe_j$, $n_{kj}$ the number of arrows in $Q$ on the shortest path from $j$ to $k$ of length $d(k,j)$, $x$ contains $n_{k,j}+\frac{l-d(k,j)}{2}$ arrows in $Q$. 

So we obtain the formula
\[
(H_A^\eta)_{k,j}=(-1)^n_{k,j}\left.\frac{H_A(t)_{k,j}}{t^{d(k,j)}}\right|_{t=\sqrt{(-1)}}
\]
and compute
\[
 H_A^\eta=
 \left(
 \begin{array}{ccccccc}
  3&0&3&0&0&0&-3\\
  0&0&0&0&0&0&0\\
  3&0&3&0&0&0&-3\\
  0&0&0&0&0&0&0\\
  0&0&0&0&0&0&0\\
  0&0&0&0&0&0&0\\
  -3&0&-3&0&0&0&3\\
 \end{array}
 \right),
\]
and its kernel is
\[
 \langle e_1+e_7,e_2,e_3+e_7,e_4,e_5,e_6\rangle.
\]

\subsection{Result}\label{HH5 results}

Now we give explicit bases for each quiver where $\psi_i\in U^*[-2]$ satisfy the properties given in section \ref{Ustar} and $\varepsilon_i\in Y^*[-h-2]$ are taken from \ref{Ystar}. 

Note the duality $HH^6(A)=(HH^5(A))^*$ which was established in \cite{EE2}, $\phi_0(z_0)\in U[-2h-2]$, $\varphi_0(\omega_i)\in Y[-h-2]$. We choose $\psi_0$ such that $\psi_0(\varphi_0(z_0))=1$ (from that follows $\psi_k(\varphi_0(z_k))=z_k\psi_k(\varphi_0(z_0))=\psi_0(\varphi_0(z_0))=1$ and $\varepsilon_i$ such that $\varepsilon_i(\phi_0(\omega_j))=\delta_{ij}$.

\subsubsection{$Q=D_{n+1}$, $n$ odd}
We define
\[
 \psi_{4k}=[(a_{n-1}^*a_{n-1}a_n^*a_n)^{\frac{n-1}{2}-k}],
\]
\[
 \varepsilon_{2k-1}=[e_{2k-1}-e_1],\, \varepsilon_{2k}=[e_{2k}], \varepsilon_n=[(e_n+e_{n+1})-e_1],k\leq \frac{n-1}{2}.
\]
\subsubsection{$Q=D_{n+1}$, $n$ even}
We define
\[
 \psi_{4k}=[a_{n-1}^*a_{n-1}(a_n^*a_na_{n-1}^*a_{n-1})^{\frac{n-2}{2}-k}],
\]
\[
\varepsilon_{2k+1}=[e_{2k+1}-e_1],\varepsilon_{2k}=[e_{2k}],1\leq k\leq \frac{n}{2}-1.
\]
\subsubsection{$Q=E_6$}
We define
\begin{eqnarray*}
 \psi_0&=&[a_3^*a_3(a_2^*a_2a_3^*a_3)^2],\\
 \psi_6&=&[-a_3^*a_3a_2^*a_2],\\
 \psi_8&=&[a_3^*a_3-a_2^*a_2],
\end{eqnarray*}
\[
 \varepsilon_3=[e_3],\,\varepsilon_6=[e_6].
\]

\subsubsection{$Q=E_7$}
We define
\begin{eqnarray*}
 \psi_0&=&[(a_4^*a_4a_3^*a_3)^4],\\
 \psi_8&=&[(a_4^*a_4a_3^*a_3)^2],\\
 \psi_{12}&=&[a_4^*a_4a_3^*a_3],
\end{eqnarray*}
\[
 \varepsilon_1=[e_1+e_7],\varepsilon_2=[e_2],\varepsilon_3=[e_3+e_7],\varepsilon_4=[e_4],\varepsilon_5=[e_5],\varepsilon_6=[e_6].
\]
\subsubsection{$Q=E_8$}
We define
\begin{eqnarray*}
 \psi_0&=&[(a_4^*a_4a_3^*a_3)^7],\\
 \psi_{12}&=&[(a_4^*a_4a_3^*a_3)^4],\\
 \psi_{20}&=&[(a_4^*a_4a_3^*a_3)^2],\\
 \psi_{24}&=&[a_4^*a_4a_3^*a_3],
\end{eqnarray*}
\[
 \varepsilon_1=[e_1],\varepsilon_2=[e_2],\varepsilon_3=[e_3],\varepsilon_4=[e_4],\varepsilon_5=[e_5],\varepsilon_6=[e_6],\varepsilon_7=[e_7],\varepsilon_8=[e_8].
\]

\section{$HH^6(A)$}
$HH^6(A)=U[-2h-2]\oplus Y[-h-2]=HH^0(A)/Im(d_6^*)$, and $Im(d_6^*)$ is spanned by the columns of the matrices $H_A^\eta$ which were computed in the previous section.

This gives us the following result:

\begin{proposition}
$HH^6(A)$ is a quotient of $HH^0(A)$. Im particular,
\[
HH^6(A)=\left\{
\begin{array}{ll} 
HH^0(A) & Q=E_6,\,E_8\\
HH^0(A)/(\sum\limits_{{i=1\atop odd}}^{n-2}\omega_i=0, \omega_n=\omega_{n+1}) & Q=D_{n+1},\,n\,\mbox{odd}\\
HH^0(A)/(\sum\limits_{{i=1\atop odd}}^{n-1}\omega_i=0), & Q=D_{n+1},\,n\,\mbox{even}\\
HH^0(A)/(\omega_1+\omega_3-\omega_7=0) & Q=E_7
\end{array}
\right..
\]
\end{proposition}

\section{Products involving $HH^0(A)=Z$}
Recall the decomposition $HH_0(A)=\mathbb{C}\oplus (U[-2])_+\oplus L[h-2]$. It is clear that the $\mathbb{C}$-part acts on $HH^i(A)$ as the usual multiplication with $\mathbb{C}$, with $z_0$ as identity. From the periodicity of the Schofield resolution with period $6$, it follows that the multiplication with $\varphi(z_0)\in HH^6(A)$ gives the natural isomorphism $HH^i(A)\rightarrow HH^{i+6}(A)$ for $i\geq1$. 

We summarize all products not involving the $\mathbb{C}$-part.
\subsection{$HH^0(A)\times HH^0(A)\twoheadrightarrow HH^0(A)$}
This is already done in the $HH^0(A)$-section of this paper. We state the results:
\subsubsection{$Q=D_{n+1}$, $n$ odd}
The products are
\[
 z_{4j}z_{4k}=\left\{\begin{array}{ll}z_{4(j+k)}&j+k<\frac{n-1}{2}\\\omega_n-\omega_{n+1}&j+k=\frac{n-1}{2}\\0&j+k>\frac{n-1}{2}\end{array}\right..
\]
\subsubsection{$Q=D_{n+1}$, $n$ even}
The products are
\[
 z_{4j}z_{4k}=\left\{\begin{array}{ll}z_{4(j+k)}&j+k<\frac{n-1}{2}\\0&j+k\geq\frac{n-1}{2}\end{array}\right..
\]
\subsubsection{$E_6$} All products are zero.
\subsubsection{$E_7$} The only nonzero product is $z_8^2=\omega_1+\omega_3-\omega_7$.
\subsubsection{$E_8$} The only nonzero product is $z_{12}^2=z_{24}$.

\subsection{$HH^0(A)\times HH^1(A)\twoheadrightarrow HH^1(A)$}
From the definition of the maps $\theta_k$ (which are generated by the central elements $z_k$), it follows that the $Z$-action is natural, i.e. the multiplication rule is the same as with the $z_k$ counterpart: $z_k\theta_0=\theta_k$. 

We state the other nonzero products:

\subsubsection{$Q=D_{n+1}$}
We have 
$z_{4j}\theta_{4k}=\theta_{4(j+k)}$ if $j+k<\frac{n-1}{2}$.

\subsubsection{$E_8$} We have $z_{12}\theta_{12}=\theta_{24}$.

\subsection{$HH^0(A)\times HH^i(A)\twoheadrightarrow HH^i(A)$, $i=2$ or $3$}
$HH^2(A)=K[-2]$ and $HH^3(A)=K^*[-2]$ live in only one degree, so $(U[-2])_+\ HH^0(A)$ acts by zero.

\subsection{$HH^0(A)\times HH^4(A)\twoheadrightarrow HH^4(A)$}
We defined $\zeta_k$, such that $z_k\zeta_k=\zeta_0$ holds. By degree arguments, only these other products are nonzero:
\subsubsection{$Q=D_{n+1}$} For $l<k$, $z_{4l}\zeta_{4k}=\zeta_{4(k-l)}$ (since $z_{4(k-l)}(z_{4l}\zeta_{4k})=(z_{4(k-l)}z_{4l})\zeta_{4k})=\zeta_0$, and $\zeta_{4(k-l)}$ is only one element of degree $-4-4(k-l)$ in $HH^4(A)$).
\subsubsection{$Q=E_8$} We have $z_{12}\zeta_{24}=\zeta_{12}$ (since $z_{12}(z_{12}\zeta_{24})=(z_{12}z_{12})\zeta_{24}=\zeta_0$, and $\zeta_{12}$ is the only element of degree $-16$ in $HH^4(A)$).

\subsection{$HH^0(A)\times HH^5(A)\twoheadrightarrow HH^5(A)$}
By definition, $z_k\psi_k=\psi_0$ holds. Since $\psi_i\in U^*[-2]$ corresponds to $\zeta_i\in U^*[-2]$ in $HH^4(A)$ with the rule $z_k\psi_k=\psi_0$ corresponding to $z_k\zeta_k=\zeta_0$ above, the multiplication rules of $\psi_k$ with elements in $HH^0(A)$ can be derived from above.

Products involving $\omega_i\in L[h-2]\subset HH^0(A)$ and\\ $\varepsilon_j=\sum\limits_{k\in F} \lambda_ke_k\in Y^*[-h-2]$ are easy to calculate: $\omega_i\varepsilon_j=\lambda_i[\omega_i]=\lambda_i\psi_0$.

\begin{proposition} 
 The multiplication $((U[-2])_+)\times Y^*[-h-2]\rightarrow HH^5(A)$ is zero
\end{proposition}
 We will show this for any quiver separately.

\subsubsection{$Q=D_{n+1}$, $n$ odd} 
For $l<k$, $z_{4l}\psi_{4k}=\psi_{4(k-l)}$. 

The nonzero products involving $\omega_i\in L[h-2]\subset HH^0(A)$ and $\varepsilon_j\in Y^*[-h-2]$ are
\[
 \omega_{2k-1}\varepsilon_{2k-1}=\omega_{2k}\varepsilon_{2k}=\omega_n\varepsilon_n=
 \omega_{n+1}\varepsilon_n=\omega_1\varepsilon_{2k-1}=\omega_1\varepsilon_n=\psi_0,
\]
\[
 \omega_1\varepsilon_{2k-1}=\omega_1\varepsilon_{n}=-\psi_0.
\]
We show $(U[-2])_+\times Y^*[-h-2]\stackrel{0}{\rightarrow} HH^5(A)$:
by degree argument, $z_{4k}\varepsilon_i=\lambda\psi_{2n-2-4k}$. Then $z_{2n-2-4k}(z_{4k}\varepsilon_i)=\lambda z_{2n-2-4k}\psi_{2n-2-4k}=\lambda\psi_0$, and  by associativity this equals 
$(z_{2n-2-4k}z_{4k})\varepsilon_i=(\omega_n-\omega_{n+1})\varepsilon_i=0$, so $\lambda=0$.

\subsubsection{$Q=D_{n+1}$, $n$ even} 
For $l<k$, $z_{4l}\psi_{4k}=\psi_{4(k-l)}$. 

The nonzero products involving $\omega_i\in L[h-2]\subset HH^0(A)$ and $\varepsilon_j\in Y^*[-h-2]$ are
\[
 \omega_{2k+1}\varepsilon_{2k+1}=\omega_{2k}\varepsilon_{2k}=\psi_0,
 \]
 \[
 \omega_1\varepsilon_{2k+1}=-\psi_0.\\
\]
We show $(U[-2])_+\times Y^*[-h-2]\stackrel{0}{\rightarrow}HH^5(A)$:
by degree argument, $z_{4k}\varepsilon_i=\lambda\psi_{2n-2-4k}$. Then 
$z_{2n-2-4k}(z_{4k}\varepsilon_i)=\lambda z_{2n-2-4k}\psi_{2n-2-4k}=\lambda\psi_0$, and this equals 
$(z_{2n-2-4k}z_{4k})\varepsilon_i=0$, so $\lambda=0$.

\subsubsection{$Q=E_6$}
The nonzero products involving $\omega_i\in L[h-2]\subset HH^0(A)$ and $\varepsilon_j\in Y^*[-h-2]$ are 
\[
 \omega_3\varepsilon_3=\omega_6\varepsilon_6=\psi_0.
\]
By degree argument, $(U[-2])_+\times Y^*[-h-2]\stackrel{0}{\rightarrow}HH^5(A)$.

\subsubsection{$Q=E_7$}
The nonzero products involving $\omega_i\in L[h-2]\subset HH^0(A)$ and $\varepsilon_j\in Y^*[-h-2]$ are 
\[
 \omega_1\varepsilon_1=\omega_2\varepsilon_2=\omega_3\varepsilon_3=\omega_4\varepsilon_4=\omega_5\varepsilon_5=\omega_6\varepsilon_6=\omega_7\varepsilon_1=\omega_7\varepsilon_3=\psi_0.
\]
We show $(U[-2])_+\times Y^*[-h-2]\stackrel{0}{\rightarrow}HH^5(A)$: by degree argument, only products involving $z_8$ may eventually be nontrivial,
\[
 z_8\varepsilon_i=\lambda \psi_8,\quad\lambda\in\mathbb{C}.
\]
Then
\[
 z_8(z_8\varepsilon_i)=\lambda z_8\psi_8=\lambda\psi_0,
\]
and by associativity this equals
\[
 z_8^2\varepsilon_i=(\omega_1+\omega_3-\omega_7)\varepsilon_i=0,
\]
so $\lambda=0$.

\section{Products involving $HH^1(A)$}

\subsection{$HH^1(A)\times HH^1(A)\stackrel{0}{\rightarrow} HH^2(A)$} This follows by degree argument since $\deg HH^1(A)>0$, $\deg HH^2(A)=-2$.

\subsection{$HH^1(A)\times HH^2(A)\twoheadrightarrow HH^3(A)$}\label{H1xH2}
$HH^2(A)$ and $HH^3(A)$ are trivial for $Q=D_{n+1}$ where $n$ is odd and for $Q=E_7,\,E_8$.

We know that $HH^1(A)$ is generated by maps $\theta_k$ and $HH^2(A)$ by $f_i$ ($i\neq\nu(i)$, and we lift 
 \begin{eqnarray*}
  f_i:A\otimes A[2]&\longrightarrow& A,\\
  1\otimes1&\longmapsto&e_i-e_{\nu(i)}
 \end{eqnarray*}
to
 \begin{eqnarray*}
  \hat f_i:A\otimes A[2]&\longrightarrow& A\otimes A,\\
  1\otimes1&\longmapsto&e_i\otimes e_i-e_{\nu(i)}\otimes e_{\nu(i)}.
 \end{eqnarray*}
 Then 
\[
 \hat f_id_3(1\otimes 1)=\hat f_i(\sum\limits_{x_j\in B} x_j\otimes x_j^*)=\sum\limits_{x_j\in B} x_je_i\otimes e_ix_j^*-x_je_{\nu(i)}\otimes e_{\nu(i)}x_j^*.
\]
To compute the lift $\Omega f_i$, we need to find out the preimage of $\sum x_je_i\otimes e_ix_j^*-x_je_{\nu(i)}\otimes e_{\nu(i)}x_j^*$ under $d_1$.

\begin{definition}
 Let $b_1,\ldots, b_k$ be arrows, $p$ the monomial $\pm b_1\cdots b_k$ and define 
 \[
v_p:=\pm (1\otimes b_1\otimes b_2\cdots b_k+b_1\otimes b_2\otimes b_3\cdots b_k+\ldots+b_1\cdots b_{k-1}\otimes b_k\otimes 1),
 \]
 and for $i<j$,
 \[
  v_p^{(i,j)}:=\pm\sum\limits_{l=i}^j b_1\cdots b_{l-1}\otimes b_l\otimes b_{l+1}\cdots b_k
 \]
\end{definition}
We will use the following lemma in our computations.
\begin{lemma}
In the above setting,
 \[
  d_1(v_p)=\pm(b_1\cdots b_k\otimes 1-1\otimes b_1\cdots b_k).
 \]
\end{lemma}
From that, we see immediately that when assuming all $x_j$ are monomials (what we can do), then
\begin{align*}
 \hat f_i(\sum\limits_{x_j\in B} x_j\otimes x_j^*)&=d_1(\sum\limits_{x_j\in B} v_{x_je_ix_j^*}^{(1,\deg(x_j))}-v_{x_je_{\nu(i)}x_j^*}^{(1,\deg(x_j))})+1\otimes\underbrace{\sum\limits_{x_j\in B} (x_je_ix_j^*-x_je_{\nu(i)}x_j^*)}_{=0},
\end{align*}
so we have
\begin{eqnarray*}
 \Omega f_i:\Omega^3(A)&\longrightarrow&\Omega(A),\\
 1\otimes 1&\longmapsto&\sum\limits_{x_j\in B} v_{x_je_ix_j^*}^{(1,\deg(x_j))}-v_{x_je_{\nu(i)}x_j^*}^{(1,\deg(x_j))}.
\end{eqnarray*}

Then 
\[
 \theta_k(\sum\limits_{x_j\in B} v_{x_je_ix_j^*}^{(1,\deg(x_j))}-v_{x_je_{\nu(i)}x_j^*}^{(1,\deg(x_j))})=z_k(\sum\limits_{x_j\in B_{-,i}}s(x_j)x_jx_j^*-\sum\limits_{x_j\in B_{-,\nu(i)}}s(x_j)x_jx_j^*),
\]
where for $s(x_j)$ is the number of arrows in $Q^*$ in the monomial expression of $x_j$.

So we get

\[
 (\theta_k\circ\Omega f_i)(1\otimes1)=z_k(\sum\limits_{x_j\in B_{-,i}}s(x_j)x_jx_j^*-\sum\limits_{x_j\in B_{-,\nu(i)}}s(x_j)x_jx_j^*).
\]

Under our identification in \cite[Subsection 4.5.]{EE2},

\[
\theta_kf_i=
[z_k(\sum\limits_{l\in I}\sum\limits_{x_j\in B_{l,i}}s(x_j)\omega_l-\sum\limits_{l\in I}\sum\limits_{x_j\in B_{l,\nu(i)}}s(x_j)\omega_l)]\in\, HH^3(A).
\]

All products are zero if $z_k$ lies in a positive degree, so we only have to calculate the products where $k=0$.

We make the following 
\begin{proposition}\label{alpha}
The multiplication with $\theta_0$ induces a symmetric isomorphism
 \[
 \alpha: HH^2(A)=K[-2]\stackrel{\cong}{\rightarrow}K^*[-2]=HH^3(A).
 \]
\end{proposition}

Now we have to work with explicit basis elements $x_j\in Ae_i$, $i\neq\nu(i)$, so we treat the Dynkin quivers separately and find the matrix $M_\alpha$ which represents this map.

\subsubsection{$Q=D_{n+1}$, $n$ even}
We can work with the basis given in section \ref{basis D_{n+1}} and compute

\begin{equation}\label{H1xH2D}
 \theta_0f_n=\frac{n}{2}([\omega_{n+1}]-[\omega_n])=-nh_n
\end{equation}
because of the relation $[\omega_n]+[\omega_{n+1}]=0$ in $HH^3(A)$. $\alpha$ is given by the matrix \[M_\alpha=(-n).\]

\subsubsection{$E_6$} 
We will write out the basis elements of $Ae_1,Ae_5$:
\begin{align*}
 B_{1,1}&=\langle e_1,a_1a_2a_5^*a_5a_2^*a_1^*\rangle,\\
 B_{2,1}&=\langle a_1^*, a_2a_5^*a_5a_2^*a_1^*, a_2a_3^*a_3a_5^*a_5a_2^*a_1^*\rangle,\\
 B_{3,1}&=\langle a_2^*a_1^*,a_3^*a_3a_2^*a_1^*,a_3^*a_3a_3^*a_3a_2^*a_1^*,a_5^*a_5a_3^*a_3a_3^*a_3a_2^*a_1^*\rangle,\\
 B_{4,1}&=\langle a_3a_2^*a_1^*,a_3a_5^*a_5a_2^*a_1^*,a_3a_5^*a_5a_3^*a_3a_5^*a_5a_2^*a_1^*\rangle,\\
 B_{5,1}&=\langle a_4a_3a_2^*a_1^*, a_4a_3 a_5^*a_5a_3^*a_3a_5^*a_5 a_2^*a_1^*\rangle,\\
 B_{6,1}&=\langle a_5a_2^*a_1^*,a_5 a_3^*a_3a_5^*a_5 a_2^*a_1^*\rangle,
\end{align*}
 and
 \[
  e_iAe_5=\langle \eta(x)|x\in e_{\nu(i)}Ae_1\rangle,
 \]
where $\eta(a)=-\epsilon_a\bar a$ and for any arrow $a:i\rightarrow j$, $\bar a$ is the arrow $j:\rightarrow i$,
so $\eta$ preserves the number of star letters of a monomial $x$. From this, we obtain

\[
 \theta_0f_1=-4[\omega_1]-2[\omega_2]+2[\omega_4]+4[\omega_5]=-8h_1-4h_2.
\]
because of the relations $\omega_1+\omega_4=\omega_2+\omega_3=0$ in $HH^3(A)$.

We do the same thing for $Ae_2$ and $Ae_4$:

\begin{align*}
 B_{1,2}&=\langle a_1, a_1a_2a_5^*a_5a_2^*,a_1a_2a_5^*a_5a_3^*a_3a_2^*\rangle,\\
 B_{2,2}&=\langle e_2, a_2a_2^*,a_2a_5^*a_5a_2^*,a_2a_3^*a_3a_5^*a_5a_2^*,a_2a_5^*a_5a_3^*a_3a_2^*, a_2a_5^*a_5a_3^*a_3a_5^*a_5a_2^*\rangle,\\
 B_{3,2}&=\langle a_2^*, a_5^*a_5a_2^*, a_3^*a_3a_2^*, a_5^*a_5a_3^*a_3a_2^*, a_3^*a_3a_5^*a_5a_2^*, \\ 
 &\quad a_3^*a_3a_5^*a_5a_3^*a_3a_2^*, a_5^*a_5a_3^*a_3a_5^*a_5a_2^*, a_5^*a_5a_3^*a_3a_5^*a_5a_3^*a_3a_2^*\rangle,\\
 B_{4,2}&=\langle a_3a_2^*, a_3a_5^*a_5a_2^*, a_3a_3^*a_3a_2^*, a_3a_3^*a_3a_5^*a_5a_2^*, \\
 &\quad a_3a_5^*a_5a_3^*a_3a_5^*a_5a_2^*, a_3a_3^*a_3a_5^*a_5a_3^*a_3a_5^*a_5a_2^*\rangle,\\
 B_{5,2}&=\langle a_4a_3a_2^*, a_4a_3a_5^*a_5a_2^*, a_4a_3a_5^*a_5a_3^*a_3a_5^*a_5a_2^*\rangle,\\
 B_{6,2}&=\langle a_5a_2^*, a_5a_3^*a_3a_2^*, a_5a_3^*a_3a_5^*a_5a_2^*, a_5a_3^*a_3a_5^*a_5a_3^*a_3a_2^*\rangle,
\end{align*}
and we get the basis elements for $e_iAe_3$ from $\eta(x_j)$ where $x_j\in e_{\nu(i)}Ae_4$. Since $\eta$ preserves the number of star-letters of a monomial, we can immediately calculate

\[
 \theta_0f_2=-2[\omega_1]-4[\omega_2]+4[\omega_4]+2[\omega_5]=-4h_1-8h_2
\]
because of the relations $[\omega_1]+[\omega_4]=[\omega_2]+[\omega_3]=0$ in $HH^3(A)$.

So $\alpha$ is given by the symmetric, nondegenerate matrix

\begin{equation}\label{H1xH2E}
 M_\alpha=\left(\begin{array}{ll}-8&-4\\-4&-8\end{array}\right).
\end{equation}

\subsection{$HH^1(A)\times HH^3(A)\stackrel{0}{\rightarrow} HH^4(A)$} This follows by degree argument: $\deg HH^1(A)\geq 0$, $\deg HH^3(A)=-2$, but $\deg HH^4(A)\leq -4$.

\subsection{$HH^1(A)\times HH^4(A)\rightarrow HH^5(A)$}
\begin{proposition}
Given $\theta_k\in HH^1(A)$ and $\zeta_l\in HH^4(A)$, we get the following cup product:
\begin{equation}\label{thetaxzeta}
 \theta_k\zeta_l=\psi_lz_k.
\end{equation}
\end{proposition}
\begin{proof}
It is enough to show $\theta_0\zeta_0=\psi_0$: $z_l(\theta_0\zeta_l)=\theta_0\zeta_0\psi_0$ implies that $(\theta_0\zeta_l)=\psi_l$, and the equation above follows from $\theta_k=z_k\theta_0$.

Let in general $x=\sum\limits_{a\in\bar Q}a\otimes x_a\in HH^4(A)$. Then $x$ represents the map
\begin{eqnarray*}
 x:=A\otimes V\otimes\mathfrak{N}h&\longrightarrow&A,\\
 1\otimes a_i\otimes 1&\longmapsto&-x_{a_i^*}\\
 1\otimes a_i^*\otimes 1&\longmapsto&x_{a_i},
\end{eqnarray*}

and it lifts to 
\begin{eqnarray*}
\hat x:A\otimes V\otimes\mathfrak{N}[h]&\longrightarrow&A\otimes A,\\
 1\otimes a\otimes 1&\longmapsto&-1\otimes x_{a^*}\\
 1\otimes a^*\otimes 1&\longmapsto&1\otimes x_{a}.
\end{eqnarray*}

Then 

\begin{eqnarray*}
 (\hat x\circ d_5)(1\otimes1)&=&\hat x(\sum\limits_{a\in\bar Q}\epsilon_a a\otimes a^*\otimes 1+\sum\limits_{a\in\bar Q}\epsilon_a 1\otimes a\otimes a^*)\\
 &=&
  \sum\limits_{a\in Q} a\otimes x_{a}-\sum\limits_{a\in Q} 1\otimes x_{a}\eta(a)+\sum\limits_{a\in Q} a^*\otimes x_{a^*}-\sum\limits_{a\in Q}1\otimes x_{a^*}\eta(a^*)\\
  &=&\sum\limits_{a\in Q} a\otimes x_{a}-\sum\limits_{a\in Q} 1\otimes ax_{a}
  \sum\limits_{a\in Q} a^*\otimes x_{a^*}-\sum\limits_{a\in Q} 1\otimes a^*x_{a^*}\\
  &=&d_1(\sum\limits_{a\in Q}1\otimes a\otimes x_{a}+1\otimes a^*\otimes x_{a^*}),
\end{eqnarray*}
so we have 
\begin{eqnarray*}
 \Omega x:\Omega^5(A)&\longrightarrow&\Omega(A),\\
 1\otimes1&\longmapsto&\sum\limits_{a\in Q}1\otimes a\otimes x_{a}+1\otimes a^*\otimes x_{a^*},
\end{eqnarray*}
and this gives us

\[
 (\theta_0\circ x)(1\otimes1)=\sum\limits_{a\in Q} a^*x_{a^*},
\]
so the cup product is 
\begin{equation}\label{thetadotx}
 \theta_0\cdot x=\sum\limits_{a\in Q} a^*x_{a^*}.
\end{equation}

It can be easily checked by using explicit elements that the RHS is $\psi_0$ for $x=\zeta_0$, but we the reason here why this is true: for $x=\sum\limits a\otimes x_a=\zeta_0$, the RHS becomes
\[
 \sum\limits_{a\in Q} a^*x_{a^*}=\sum\limits_{a\in Q} (a^*,x_{a^*})[\omega_{t(a)}],
\]
where $(-,-): A\times A\rightarrow\mathbb{C}$ is the bilinear form attached to $A$ as a Frobenius algebra (see \ref{Frobenius}).

But under the bilinear form on $V\otimes A$, given in \cite[Subsection 4.3.]{EE2} which induces the duality $HH^4(A)=(HH^1(A))^*$, 
\[(a\otimes x_a,b\otimes x_b)=\delta_{a,b^*}\epsilon_a(x_a,x_b),\]
\[
 \sum\limits_{a\in Q} (a^*,x_{a^*})=(\theta_0,\zeta_0)=1.
\]
So for $x=\zeta_0$, equation (\ref{thetadotx}) becomes
\begin{equation}
\theta_0\zeta_0=(\theta_0,\zeta_0)\psi_0=\psi_0,
\end{equation}
because $[\omega_i]=\psi_0$ in $HH^5(A)$ for all $i\in I$.
\end{proof}
\subsection{$HH^1(A)\times HH^5(A)\rightarrow HH^6(A)$}\label{H1xH5}
We know that 
\begin{eqnarray*}
0&\leq&\deg(HH^1(A))\leq h-4,\\
-h-2&\leq&\deg(HH^5(A))\leq-2,\\
-2h&\leq&\deg(HH^6(A))\leq-h-2,
\end{eqnarray*}
so the product is trivial unless we pair the lowest degree parts of $HH^1(A)$ (generated by $\theta_0$) and $HH^5(A)$ (which is $Y^*[-h-2]$). The product will then live in degree $-h-2$ which is the top degree part of $HH^6(A)$, the space $Y[-h-2]$.

Given an element $\psi\in HH^5(A) (-h-2))$ which has the form
\[
\begin{array}{rcl}
 \psi:A\otimes\mathfrak{N}[h+2]&\longrightarrow&A,\\
 1\otimes1&\longmapsto& \sum\limits_{i\in F}\lambda_ie_i\in R,
\end{array}
\]
this lifts to
\[
\begin{array}{rcl}
 \hat\psi:A\otimes\mathfrak{N}[h+2]&\longrightarrow&A\otimes A,\\
 1\otimes1&\longmapsto& \sum\limits_{i\in F}\lambda_ie_i\otimes e_i.
\end{array}
\]

Then 
\begin{eqnarray*}
\hat\psi(d_6(1\otimes1))&=&\hat\psi(\sum\limits_{x_j\in B}\sum\limits_{i\in F}x_j\otimes x_j^*)=\hat\psi(\sum\limits_{x_j\in B}\sum\limits_{i\in F}\eta(x_j)\otimes\eta(x_j^*))\\
&=&\sum\limits_{x_j\in B}\sum\limits_{i\in F}\lambda_i\eta(x_j)e_i\otimes e_ix_j^*\\
&=&d_1(\sum\limits_{i\in F}\sum\limits_{x_j\in B}\lambda_i v^{(1,\deg(x_j))}_{\eta(x_j)e_ix_j^*})+1\otimes\underbrace{\sum\limits_{i\in F}\lambda_i\eta(x_j)e_ix_j^*}_{=0},
\end{eqnarray*}
so $\psi$ lifts to
\begin{eqnarray*}
 \Omega\psi:\Omega^6(A)&\longrightarrow&\Omega(A),\\
 1\otimes1&\longmapsto&\sum\limits_{i\in F}\sum\limits_{x_j\in B}\lambda_i v^{(1,\deg(x_j))}_{\eta(x_j)e_ix_j^*}.
\end{eqnarray*}

We get
 
\[
 (\theta_0\circ\Omega\psi)(1\otimes1)=\sum\limits_{i\in F}\sum\limits_{x_j\in Ae_i}\lambda_is(x_j)\eta(x_j)x_j^*,
\]
where $s(x_j)$ is the number of arrows in $Q^*$ in the monomial expression of $x_j$ (or in general if $x_j$ is a homogeneous polynomial where each monomial term has the same number of arrows in $Q^*$, then $s(x_j)$ is the number of $Q^*$-arrows in each monomial term).\\

Under our identifications in \cite[Subsection 4.5.]{EE2},

\[
\theta_0\psi=\sum\limits_{i\in F}\sum\limits_{x_j\in Ae_i}\lambda_is(x_j)\eta(x_j)x_j^*=\sum\limits_{i,k\in F}\sum\limits_{x_j\in e_kAe_i}\lambda_is(x_j)\eta(x_j)x_j^*
\]

To simplify this computation, we will choose a basis, such that all $x_j\in e_kAe_l$ for some $k,l\in I$ and that additionally $x_j$ is an eigenvector of $\eta$ for $k,l\in F$ (since $\eta$ is an involution on $e_kAe_l$ for $k,l\in F$).  Let $B_{k,l}^+$ be a basis of $(e_kAe_l)_+=\ker(\eta|_{e_kAe_l}-1)$ and $B_{k,l}^-$ a basis of $(e_kAe_l)_-=\ker(\eta|_{e_kAe_l}+1)$.

Let us define
\begin{equation}
 \kappa_{k,l}=\sum\limits_{x_j\in B_{k,l}^+}s(x_j)-\sum\limits_{x_j\in B_{k,l}^-}s(x_j).
\end{equation}

Then the above equation becomes 

\begin{equation}
\theta_0\psi=\sum\limits_{l\in F}\lambda_l\sum\limits_{k\in F}\kappa_{k,l}\varphi_0(\omega_k).
\end{equation}

\begin{proposition}
 The multiplication by $\theta_0$ induces a skew-symmetric isomorphism
 \[
 \beta: Y^*[-h-2]\stackrel{\cong}{\rightarrow} Y[-h-2].
 \]
\end{proposition}
We will treat the Dynkin quivers separately and find the matrix $M_\beta$ which represents $\beta$ for each of these quivers.

\subsubsection{$Q=D_{n+1}$, $n$ odd}
We use the same basis as given in section \ref{basis D_{n+1}}. Recall that these basis elements have the property $\eta(x)=(-1)^{n_x}x$ where $n_x$ is the number of $Q$-arrows in the monomial expression of $x$. 

We can compute that for $k\leq l\leq n-1$,
\begin{eqnarray*}
\kappa_{k,l} &=& 
\left\{\begin{array}{lc}
\sum\limits_{j=0}^{k-1}(-1)^j(l-k+j)+(-1)^{n-l}\sum\limits_{j=0}^{k-1}(-1)^j(n-k+j),&k\leq n-l\\
\begin{array}{l}
\sum\limits_{j=0}^{n-l-1}(-1)^j(l-k+j)+(-1)^{n-l}\sum\limits_{j=0}^{k-1}(-1)^j(n-k+j)\\
+(-1)^{n-l}\sum\limits_{j=0}^{k-1+l-n}(-1)^j(n-k+j),
\end{array}
&k>n-l
\end{array}
\right.\\
&=&
 \left\{
 \begin{array}{cc} 
  n-k+l-1&k\,odd,\quad l\,odd\\
  l-n&k\,odd,\quad l\,even\\
  -k&k\,even,\quad l\,odd\\
  0&k\,even,\quad l\,even
  \end{array}
  \right.
\end{eqnarray*}
for $l\leq k\leq n-1$,
\begin{eqnarray*}
\kappa_{k,l}&=& 
\left\{
\begin{array}{lc}
 (-1)^{k-l}\sum\limits_{j=0}^{l-1}(-1)^jj+(-1)^{n-l}\sum\limits_{j=0}^{l-1}(-1)^j(n-k+j)&j\leq n-k\\
 \begin{array}{l}
 (-1)^{k-l}\sum\limits_{j=0}^{n-k-1}(-1)^jj+(-1)^{n-l}\sum\limits_{j=0}^{l-1}(-1)^j(n-k+j)\\
 +(-1)^{n-l}\sum\limits_{j=0}^{l+k-n-1}(-1)^j(n-k+j)
 \end{array}
 &j>n-k
\end{array}
\right.\\
 &=&
 \left\{
 \begin{array}{cl} 
  n-k+l-1&k\,odd,\quad l\,odd\\
  l&k\,odd,\quad l\,even\\
  n-k&k\,even,\quad l\,odd\\
  0&k\,even,\quad l\,even
  \end{array}
  \right.
\end{eqnarray*}
for $k,l\leq n-1$,
\begin{eqnarray*}
  \kappa_{k,n}=\kappa_{k,n+1}&=&
  \sum\limits_{j=0}^{k-1}(-1)^j(n-k-j)=
  \left\{
  \begin{array}{cl}
   n-\frac{k+1}{2}&k\,odd\\
   -\frac{k}{2}&k\,even
  \end{array}
  \right.\\
  \kappa_{n,l}=\kappa_{n+1,l}&=&
  (-1)^{n-j}\sum\limits_{j=0}^{l-1}(-1)^jj=
  \left\{
  \begin{array}{cl}
   n-\frac{l-1}{2}&l\,odd\\
   \frac{l}{2}&k\,even\\
  \end{array}
  \right.\\
  \kappa_{n,n}=\kappa_{n+1,n+1}&=&\sum\limits_{j=0}^{\frac{n-1}{2}}2j=\frac{n^2-1}{4}\\
 \kappa_{n+1,n}=\kappa_{n,n+1}&=&\sum\limits_{j=0}^{\frac{n-3}{2}}(-1\cdot(1+2l))=-\left(\frac{n-1}{2}\right)^2
\end{eqnarray*}

$Y^*[-h-2]$ has basis $\varepsilon_{2k+1}=[e_{2k+1}-e_1]$ ($0\leq k\leq \frac{n-3}{2}$), $\varepsilon_{2k}=[e_{2k}]$ ($k\leq\frac{n-1}{2}$), $\varepsilon_n=[e_n+e_{n+1}-e_1]$, and we can calculate the products
\begin{eqnarray*}
\theta_0\varepsilon_{2k+1}&=&\sum\limits_{i\in F}(\kappa_{i,2k+1}-\kappa_{i,1})\varphi_0(\omega_i)\\&=&
2k\sum\limits_{{i=1\atop odd}}^{n-2}\varphi_0(\omega_i)-n\sum\limits_{{i=2\atop even}}^{2k}\varphi_0(\omega_i)+k\varphi_0(\omega_n+\omega_{n+1}),\\
\theta_0\varepsilon_{2k}&=&\sum\limits_{i\in F}(\kappa_{i,2k+1})\varphi_0(\omega_i)\\
&=&(2k-n)\sum\limits_{{i=1\atop odd}}^{2k-1}\varphi_0(\omega_i)+2k\sum\limits_{{i=2k+1\atop odd}}^{n-2}\varphi_0(\omega_i)+k\varphi_0(\omega_n+\omega_{n+1}),\\
\theta_0\varepsilon_n&=&\sum\limits_{i\in F}(\kappa_{i,n}+\kappa_{n+1,1}-\kappa_{i,1})\varphi_0(\omega_i)\\
&=&(n-1)\sum\limits_{i=1,odd}^{n-2}\varphi_0(\omega_i)-n\sum\limits_{{i=2\atop even}}^{n-1}\varphi_0(\omega_i)+\frac{n-1}{2}\varphi_0(\omega_n+\omega_{n+1})
\end{eqnarray*}
We use the defining relations in $Y[-h-2]$, 
\begin{eqnarray*}
 \varphi_0(\omega_1)&=&-\varphi_0(\sum\limits_{{i=3\atop odd}}^{n-2}\varphi_0(\omega_i)-\varphi_0(\omega_n))\\
\varphi_0(\omega_{n+1})&=&\varphi_0(\omega_n)
\end{eqnarray*}
to write the RHS of the above cup product calculations in terms of the basis $(\omega_i)_{2\leq i\leq n}$:
\begin{eqnarray*}
\theta_0\varepsilon_{2k+1}&=&-n\sum\limits_{{i=2\atop even}}^{2k}\varphi_0(\omega_i),\\ 
\theta_0\varepsilon_{2k}&=&n\sum\limits_{{i=2k+1\atop  odd}}^{n-2}\varphi_0(\omega_i)+n\varphi_0(\omega_n),\\
\theta_0\varepsilon_n&=&-n\sum\limits_{{i=2\atop even}}^{n-1}\varphi_0(\omega_i).
\end{eqnarray*}

$\beta$ is given by the skew-symmetric, nondegenerate matrix 
\[
 M_\beta=\left(
 \begin{array}{ccccccccccc}
  0&-n&0&-n&\ldots&\ldots&-n&0&-n&0&-n\\
  n&0&0&0&\ldots&\ldots&0&0&0&0&0\\
  0&0&0&-n&\ldots&\ldots&-n&0&-n&0&-n\\
  n&0&n&0&\ldots&\ldots&0&0&0&0&0\\
  0&0&0&0&\ldots&\ldots&-n&0&-n&0&-n\\
  \vdots&\vdots&\vdots&\vdots&\ddots&&\vdots&\vdots&\vdots&\vdots&\vdots\\
  \vdots&\vdots&\vdots&\vdots&&\ddots&\vdots&\vdots&\vdots&\vdots&\vdots\\
  n&0&n&0&\ldots&\ldots&0&0&0&0&0\\
  0&0&0&0&\ldots&\ldots&0&0&-n&0&-n\\
  n&0&n&0&\ldots&\ldots&0&n&0&0&0\\
  0&0&0&0&\ldots&\ldots&0&0&0&0&-n\\
  n&0&n&0&\ldots&\ldots&0&n&0&n&0\\      
 \end{array}
 \right)
\]
with respect to the chosen basis $\varepsilon_2,\varepsilon_3,\ldots\varepsilon_n$ of $Y^*[-h-2]$ and the dual basis $\varphi_0(\omega_2),\varphi_0(\omega_3)\ldots\varphi_0(\omega_n)$ of $Y[-h-2]$).

\subsubsection{$Q=D_{n+1}$, $n$ even}
We use the same basis as in section \ref{D_{n+1}, n even} for our computations.

For $k\leq l\leq n-1$,
\begin{eqnarray*}
\kappa_{k,l} &=& 
\left\{\begin{array}{lc}
\sum\limits_{j=0}^{k-1}(-1)^j(l-k+j)+(-1)^{n-l+1}\sum\limits_{j=0}^{k-1}(-1)^j(n-k+j),&k\leq n-l\\
\begin{array}{l}
\sum\limits_{j=0}^{n-l-1}(-1)^j(l-k+j)+(-1)^{n-l-1}\sum\limits_{j=0}^{k-1}(-1)^j(n-k+j)\\
+(-1)^{n-l}\sum\limits_{j=0}^{k-1+l-n}(-1)^j(n-k+j),
\end{array}
&k>n-l
\end{array}
\right.\\
&=&
 \left\{
 \begin{array}{cc} 
  n-k+l-1&k\,odd,\quad l\,odd\\
  l-n&k\,odd,\quad l\,even\\
  -k&k\,even,\quad l\,odd\\
  0&k\,even,\quad l\,even
  \end{array}
  \right.
\end{eqnarray*}
for $l\leq k\leq n-1$,
\begin{eqnarray*}
\kappa_{k,l}&=& 
\left\{
\begin{array}{lc}
 (-1)^{k-l}\sum\limits_{j=0}^{l-1}(-1)^jj+(-1)^{n-l+1}\sum\limits_{j=0}^{l-1}(-1)^j(n-k+j)&j\leq n-k\\
 \begin{array}{l}
 (-1)^{k-l}\sum\limits_{j=0}^{n-k-1}(-1)^jj+(-1)^{n-l+1}\sum\limits_{j=0}^{l-1}(-1)^j(n-k+j)\\
 +(-1)^{n-l}\sum\limits_{j=0}^{l+k-n-1}(-1)^j(n-k+j)
 \end{array}
 &j>n-k
\end{array}
\right.\\
 &=&
 \left\{
 \begin{array}{cl} 
  n-k+l-1&k\,odd,\quad l\,odd\\
  l&k\,odd,\quad l\,even\\
  n-k&k\,even,\quad l\,odd\\
  0&k\,even,\quad l\,even
  \end{array}
  \right.
\end{eqnarray*}

$Y^*[-h-2]$ has basis $\varepsilon_{2k}=[e_{2k}]$, $\varepsilon_{2k+1}=[e_{2k+1}-e_1]$ ($1\leq k\leq\frac{n-2}{2}$), and we calculate the products

\begin{eqnarray*}
 \theta_0\varepsilon_{2k+1}&=&\sum\limits_{i\in F}(\kappa_{i,2k+1}-\kappa_{i,1})\varphi_0(\omega_i)\\
 &=&2k\sum\limits_{{i=1,\atop odd}}^{n-1}\varphi_0(\omega_i)-n\sum\limits_{{i=2\atop even}}^{2k}\varphi_0(\omega_i),\\
\theta_0\varepsilon_{2k}&=&\sum\limits_{i\in F}(\kappa_{i,2k})\varphi_0(\omega_i)\\
 &=&(2k-n)\sum\limits_{{i=1\atop odd}}^{2k-1}[\omega_i]+2k\sum\limits_{i=2k+1}^{n-2}\varphi_0(\omega_i),
\end{eqnarray*}
and we use the defining relation of $Y[-h-2]$,
\[
 \varphi_0(\omega_1)=-\sum\limits_{{i=3\atop odd}}^{n-2}\varphi(\omega_i)
\]
to write the results of the cup product calculations in terms of the basis $\varphi_0(\omega_2),\varphi_0(\omega_3),\ldots,\varphi_0(\omega_{n-1})$. We get

\begin{eqnarray*}
 \theta_0\varepsilon_{2k+1}&=&-n\sum\limits_{{i=2\atop even}}^{2k}\varphi_0(\omega_i)\\
 \theta_0\varepsilon_{2k}&=&n\sum\limits_{{i=2k+1\atop odd}}^{n-1}\varphi_0(\omega_i).
\end{eqnarray*}

$\beta$ is given by the matrix

\[
 M_\beta=\left(
 \begin{array}{ccccccccccc}
  0&-n&0&-n&\ldots&\ldots&-n&0&-n&0&-n\\
  n&0&0&0&\ldots&\ldots&0&0&0&0&0\\
  0&0&0&-n&\ldots&\ldots&-n&0&-n&0&-n\\
  n&0&n&0&\ldots&\ldots&0&0&0&0&0\\
  0&0&0&0&\ldots&\ldots&-n&0&-n&0&-n\\
  \vdots&\vdots&\vdots&\vdots&\ddots&&\vdots&\vdots&\vdots&\vdots&\vdots\\
  \vdots&\vdots&\vdots&\vdots&&\ddots&\vdots&\vdots&\vdots&\vdots&\vdots\\
  n&0&n&0&\ldots&\ldots&0&0&0&0&0\\
  0&0&0&0&\ldots&\ldots&0&0&-n&0&-n\\
  n&0&n&0&\ldots&\ldots&0&n&0&0&0\\
  0&0&0&0&\ldots&\ldots&0&0&0&0&-n\\
  n&0&n&0&\ldots&\ldots&0&n&0&n&0\\      
 \end{array}
 \right)
\]
with respect to the basis $\varepsilon_2,\varepsilon_3,\ldots,\varepsilon_{n-1}$ and its dual basis $\varphi_0(\omega_2),\varphi_0(\omega_3),\ldots\varphi_0(\omega_{n-1})$.

\subsubsection{$Q=E_6$}
We work with the bases 
\begin{eqnarray*}
 B_{3,3}^+&=&\{ e_3,a_3^*a_3-a_2^*a_2,(a_3^*a_3-a_2^*a_2)^2,a_5^*a_5a_3^*a_3a_5^*a_5,\\
 &&a_5^*a_5a_3^*a_3a_5^*a_5a_3^*a_3,a_3^*a_3a_5^*a_5a_3^*a_3a_5^*a_5a_3^*a_3\},\\
 B_{3,3}^-&=&\{ a_5^*a_5,a_3^*a_3a_5^*a_5,a_5^*a_5a_3^*a_3,a_3^*a_3a_5^*a_5a_3^*a_3,\\
 &&a_5^*a_5a_3a_3^*(a_3^*a_3-a_2^*a_2)^2,a_3^*a_3a_5^*a_5a_3a_3^*(a_3^*a_3-a_2^*a_2)^2\},\\
 B_{6,3}^+&=&\{a_5a_3^*a_3a_5^*a_5,a_5a_3^*a_3a_5^*a_5a_3^*a_3),\\
 &&a_5a_3^*a_3a_5^*a_5a_3a_3^*(a_3^*a_3-a_2^*a_2)\},\\
 B_{6,3}^-&=&\{a_5,a_5a_3^*a_3,a_5a_3^*a_3(a_3^*a_3-a_2^*a_2)\},\\ 
 B_{3,6}^+&=&=\{a_5^*,a_3^*a_3a_5^*,(a_3^*a_3-a_2^*a_2)a_3^*a_3a_5^*\},\\
 B_{3,6}^-&=&\{a_5^*a_5a_3^*a_3a_5^*,a_3^*a_3a_5^*a_5a_3^*a_3a_5^*,a_3^*a_3a_5^*a_5(a_3^*a_3)^2a_5^*\},\\
 B_{6,6}^+&=&\{e_6,a_5a_3^*a_3a_5^*a_5(a_3^*a_3)^2a_5^*\},\\
 B_{6,6}^-&=&\{a_5a_3^*a_3a_5^*,a_5(a_3^*a_3)^2a_5^*\}.
\end{eqnarray*}
We immediately get the matrix
\[
M_\beta=\left(
\begin{array}{cc}
\kappa_{3,3}&\kappa_{3,6}\\
\kappa_{6,3}&\kappa_{6,6}
\end{array}
\right)
=
\left(
\begin{array}{cc}
0&-6\\
6&0 
\end{array}
\right)
\]
which represents the $\beta$ with respect to the basis $\varepsilon_3, \varepsilon_6$ and dual basis $\varphi_0(\omega_3),\varphi_0(\omega_6)$.

\subsubsection{$E_7$}
For $E_7$ and $E_8$ we don't have to work with an explicit basis to calculate $\kappa_{k,l}$ since for any basis element $x$, $\eta(x)=\pm x$. It is enough to know the following:
given any monomial $x\in e_kAe_j$ of length $l$, $n_{k,j}$ the number of arrows $x\in Q$ and $d(k,j)$ the distance between the vertices $k,j$, we know that 
$x$ contains $n_{k,j}+\frac{l-d(k,j)}{2}$ arrows in $Q$ and  $d(k,j)-n_{k,j}+\frac{l-d(k,j)}{2}$ arrows in $\bar Q$.

We can derive the following formula:

\begin{equation}\label{kappa E7}
 \kappa_{k,j}=(-1)^{n_{k,j}}\left(\left.(d(k,j)-n_{k,j})\frac{H_A(t)}{t^{d(k,j)}}\right|_{t=\sqrt{-1}}+\left.\frac{1}{2}t\frac{d}{dt}\frac{H_A(t)}{t^{d(k,j)}}\right|_{t=\sqrt{-1}}\right).
\end{equation}
The resulting matrix is
\[
(\kappa_{k,j})_{k,j}=
\left(
  \begin{array}{ccccccc}
  12&6&9&3&0&3&-9\\
  -6&0&3&0&0&0&-3\\
  15&-3&12&3&0&3&-12\\
  -3&0&-3&0&0&0&-6\\
  0&0&0&0&0&-9&0\\
  -3&0&-3&0&9&0&-6\\
  -15&3&-12&6&0&6&12                       
  \end{array}
\right)
\]
A basis of $Y^*[-h-2]$ is given by \[\varepsilon_1=[e_1+e_7],\varepsilon_2=[e_2],\varepsilon_3=[e_3+e_7],\varepsilon_4=[e_4],\varepsilon_5=[e_5],\varepsilon_6=[e_6],\]

$(\theta_0\varepsilon_i)_{1\leq1\leq6}$ is given by

\[
\left(
  \begin{array}{cccccc}
  3&6&0&3&0&3\\
  -9&0&0&0&0&0\\
  3&-3&0&3&0&3\\
  -9&0&-9&0&0&0\\
  0&0&0&0&0&-9\\
  -9&0&-9&0&9&0\\
  -3&3&0&6&0&6                       
  \end{array}
\right)
\left(
\begin{array}{cccccc}
 \varphi_0(\omega_1)\\
 \varphi_0(\omega_2)\\
 \varphi_0(\omega_3)\\
 \varphi_0(\omega_4)\\
 \varphi_0(\omega_5)\\
 \varphi_0(\omega_6)\\
 \varphi_0(\omega_7)  
\end{array}
\right)
\]

Now use the defining relation of $Y[-h-2]$,
\[
  \varphi_0(\omega_7)=\varphi_0(\omega_1)+\varphi_0(\omega_3)
\]
to obtain the matrix
\[
M_\beta=\left(
\begin{array}{cccccc}
 0& 9& 0& 9& 0& 9\\
-9& 0& 0& 0& 0& 0\\
 0& 0& 0& 9& 0& 9\\
-9& 0&-9& 0& 0& 0\\
 0& 0& 0& 0& 0&-9\\
-9& 0&-9& 0& 9& 0  
\end{array}
\right)
\]
which represents $\beta$ with respect to the basis $\varepsilon_1\ldots\varepsilon_6$ and its dual basis $\varphi_0(\omega_1),\ldots,\varphi_0(\omega_6)$.

\subsubsection{$E_8$}
We can use \ref{kappa E7} and get the matrix

\[
 M_\beta=(\kappa_{k,j})_{k,j}=
 \left(
 \begin{array}{cccccccc}
  0& 15&  0& 15&  0&  0&  0&-15\\
-15&  0&  0&  0&  0&  0&  0&  0\\
  0&  0&  0& 15&  0&  0&  0&-15\\
-15&  0&-15&  0&  0&  0&  0&  0\\
  0&  0&  0&  0&  0&  0&  0&-15\\
  0&  0&  0&  0&  0&  0&-15&  0\\
  0&  0&  0&  0&  0& 15&  0&-15\\
 15&  0& 15&  0& 15&  0&15&  0 
 \end{array}
 \right)
\]

which represents $\beta$ with respect to the basis $\varepsilon_1,\ldots \varepsilon_8$ and its dual basis $\varphi_0(\omega_1),\ldots,\varphi_0(\omega_8)$.\\

\begin{remark}
With respect to our chosen bases $(\varepsilon_i)_{i\in I'}$ and $\phi_0(\omega_i)_{i\in I'}$, such that the vertex set $I'\subset I$, together with the arrows in $I$ form a connected subquiver $\bar{Q}'$, $M_\beta$ can be written in this general form:
\begin{equation}
 M_\beta=\frac{h}{2}\cdot (C')^\epsilon,
\end{equation}
where we call $(C')^\epsilon$ the \emph{signed adjacency matrix} of the subquiver $\bar{Q}'$, that is 
\begin{equation}(C')_{ij}=\left\{\begin{array}{cc}
0&\mbox{if $i,j$ are not adjacent},\\
+1&\mbox{if arrow $i\leftarrow j$ lies in $Q^*$},\\
-1&\mbox{if arrow $i\leftarrow j$ lies in $Q$},\\                                    \end{array}
\right.
\end{equation}
In the $D_{n+1}$-case, we have
\[
 M_\beta=n\cdot\left(
 \begin{array}{cccccccc}
  0&1&0&\cdots&\cdots&\cdots&\cdots&0\\
  -1&0&1&0&\cdots&\cdots&\cdots&0\\
  0&-1&0&1&0&\cdots&\cdots&0\\
  \vdots&\ddots&\ddots&\ddots&\ddots&\ddots&&\vdots\\
  \vdots&&\ddots&\ddots&\ddots&\ddots&\ddots&\vdots\\
  \vdots&&&\ddots&\ddots&\ddots&\ddots&0\\
  \vdots&&&&\ddots&\ddots&\ddots&1\\
  0&\cdots&\cdots&\cdots&\cdots&0&-1&0
 \end{array}
 \right)^{-1},
\]
in the $E_6$-case, we have
\[
 M_\beta=6\cdot\left(\begin{array}{cc}0&1\\-1&0\end{array}\right)^{-1},
\]
in the $E_7$-case, we have

\[
 M_\beta=9\cdot\left(\begin{array}{cccccc}
 0&-1&0&0&0&0\\
 1&0&-1&0&0&0\\
 0&1&0&-1&0&0\\
 0&0&1&0&1&0\\
 0&0&0&-1&0&1\\
 0&0&0&0&-1&0
 \end{array}\right)^{-1},
\]
and in the $E_8$-case, we have
\[
 M_{\beta}=15\cdot\left(\begin{array}{cccccccc}
0&-1&0&0&0&0&0&0\\
1&0&-1&0&0&0&0&0\\
0&1&0&-1&0&0&0&0\\
0&0&1&0&-1&0&0&0\\
0&0&0&1&0&1&0&1\\
0&0&0&0&-1&0&1&0\\
0&0&0&0&0&-1&0&0\\
0&0&0&0&-1&0&0&0
\end{array}\right)^{-1}.
\]
\end{remark}

\newpage
\section{Products involving $HH^2(A)$}
We start with $HH^2(A)\times HH^3(A)\rightarrow HH^5(A)$ first and then deduce $HH^2(A)\times HH^2(A)\rightarrow HH^4(A)$ from associativity.

\subsection{$HH^2(A)\times HH^3(A)\rightarrow HH^5(A)$}
We will prove the following general proposition:
\begin{proposition}\label{perfect}
 For the basis elements $f_i\in HH^2(A)$, $h_j\in HH^3(A)$, the cup product is
 \begin{equation}
  f_ih_j=\delta_{ij}\psi_0.
 \end{equation}
\end{proposition}
\begin{proof}
 Recall the maps 
 \begin{eqnarray*}
 h_j:A\otimes\mathfrak{N}&\rightarrow&A,\\
 1\otimes1&\mapsto&\omega_j
 \end{eqnarray*}
 
 and lift it to 
 \begin{eqnarray*}
 \hat h_j:A\otimes\mathfrak{N}&\rightarrow&A\otimes A,\\
 1\otimes1&\mapsto&1\otimes\omega_j.
 \end{eqnarray*}
 Then 
 \[
 \hat h_j(d_4(1\otimes a\otimes 1))=\hat h_j(a\otimes1-1\otimes a)=a\otimes\omega_j=d_1(1\otimes a\otimes\omega_j),
 \]
so
\begin{eqnarray*}
 \Omega h_j:\Omega^4(A)&\rightarrow&\Omega(A),\\
 1\otimes a\otimes 1&\mapsto&1\otimes a\otimes\omega_j.
\end{eqnarray*}
Then we have
\begin{eqnarray*}
 \Omega h_j(d_5(1\otimes1))&=&\Omega h_j(\sum\limits_{a\in\bar Q Q}\epsilon_aa\otimes a^*\otimes1+\sum\limits_{a\in\bar Q}\epsilon_a1\otimes a\otimes a^*)\\
 &=&\sum\limits_{a\in\bar Q}\epsilon_a\otimes a^*\otimes\omega_j=d_2(1\otimes\omega_j),
\end{eqnarray*}
so
\begin{eqnarray*}
 \Omega^2 h_j:\Omega^5(A)&\rightarrow&\Omega^2(A),\\
 1\otimes 1&\mapsto&1\otimes \omega_j.
\end{eqnarray*}
This gives us 
\[
 f_i(\Omega^2 h_j)(1\otimes1)=f_i(1\otimes\omega_j)=\delta_{ij}\omega_j,
\]
i.e. the cup product
\[
 f_ih_j=\delta_{ij}[\omega_j]=\delta_{ij}\psi_0.
\]
\end{proof}

\subsection{$HH^2(A)\times HH^2(A)\rightarrow HH^4(A)$}
Since $\deg HH^2(A)=-2$, their product has degree $-4$ (i.e. lies in $span(\zeta_0)$, so it can be written as

\begin{eqnarray*}
 HH^2(A)\times HH^2(A)&\rightarrow&HH^4(A),\\
 (a,b)&\mapsto&\langle-,-\rangle\zeta_0,
\end{eqnarray*}
where $\langle-,-\rangle:HH^2(A)\times HH^2(A)\rightarrow\mathbb{C}$ is a bilinear form. We prove the following proposition:

\begin{proposition}
 The cup product $HH^2(A)\times HH^2(A)\rightarrow HH^4(A)$ is given by $\langle-,-\rangle=\alpha$, where $\alpha$ (from Proposition \ref{alpha}) is regarded as a syymetric bilinear form.
\end{proposition}
\begin{proof}

We use (\ref{thetaxzeta}) to get
\begin{equation}\label{H2xH2 1}
 \theta_0(f_if_j)=\theta_0(\langle f_i,f_j\rangle\zeta_0)=\langle f_i,f_j\rangle\psi_0.
\end{equation}

On the other hand, by Proposition \ref{alpha} and Proposition \ref{perfect},
\begin{equation}\label{H2xH2 2}
 (\theta_0f_i)f_j=\alpha(f_i)f_j=\sum (M_\alpha)_{li}h_lf_j=(M_\alpha)_{ji}\psi_0=(M_\alpha)_{ij}\psi_0.
\end{equation}
By associativity of the cup product, we can equate (\ref{H2xH2 1}) and (\ref{H2xH2 2}) to get
\begin{equation}
 \langle f_i,f_j\rangle=(M_\alpha)_{ij}.
\end{equation}
\end{proof}

\newpage

\subsection{$HH^2(A)\times HH^4(A)\stackrel{0}{\rightarrow}HH^6(A)$}
This computation uses the Batalin-Vilkovisky structure on Hochschild cohomology: We have $\deg HH^2(A)=-2$, $\deg HH^4(A)\geq -h$ and $\deg HH^6(A)\leq -h-2$. So we know by degree argument that 
\begin{equation}
 f_k\zeta_l=\left\{\begin{array}{cc}0&l>h-4\\\sum\limits_s\lambda_s\varphi(\omega_s)&l=h-4\end{array}\right..
\end{equation}
We use \cite[(6.0.12)]{Eu3} and the isomorphism $HH^i(A)=HH_{6m+2-i}(A)$ to get for the Gerstenhaber bracket on $HH^*(A)$:

\begin{eqnarray*}
 [f_k,\zeta_l]&=&\Delta(f_k\zeta_l)-\underbrace{\Delta(f_k)}{=0}\zeta_l-f_k\underbrace{\Delta(\zeta_l)}_{=0}\\
 &=&\sum\limits_s\lambda_s(\frac{1}{2}+m)h\beta^{-1}(\varphi(\omega_s))
\end{eqnarray*}
The Gerstenhaber bracket has to be independent of the choice of $m\geq0$. This implies that the RHS has to be zero, so all $\lambda_s=0$. This shows that 
\begin{equation}
 f_k\zeta_{h-4}=0,
\end{equation}
so we have that the cup product of $HH^2(A)$ with $HH^4(A)$ is zero.

\subsection{$HH^2(A)\times HH^5(A)\stackrel{0}{\rightarrow}HH^7(A)$}
Let $a\in HH^2(A)$ and $b\in HH^5(A)$ be homogeneous elements, then $ab=\lambda \theta_k\in HH^7(A)= U[-2h-2]$, $\lambda\in\mathbb{C}$. Then 
\[
 \lambda\psi_0=\lambda\psi_kz_k=\lambda\theta_k\zeta_k=\lambda b(a\zeta_k)=0,
\]
the last equality coming from the product $a\zeta_k\in HH^2(A)\cup HH^4(A)=0$.

\section{Products involving $HH^3(A)$}
\subsection{$HH^3(A)\times HH^3(A)\stackrel{0}{\rightarrow}HH^6(A)$} This follows by degree argument: $\deg HH^3(A)=-2$, $\deg HH^6(A)\leq -h-2<-4$.
\subsection{$HH^3(A)\times HH^4(A)\stackrel{0}{\rightarrow}HH^7(A)$} This follows by degree argument: $\deg HH^3(A)=-2$, $\deg HH^4(A)\geq -h$, $\deg HH^7(A)\leq -h-4<-h-2$.
\subsection{$HH^3(A)\times HH^5(A)\stackrel{0}{\rightarrow}HH^8(A)$} This follows by degree argument: $\deg HH^3(A)=-2$, $\deg HH^5(A)\geq -h-2$, $\deg HH^8(A)=-2h-2<-h-4$.

\section{Products involving $HH^4(A)$}
\subsection{$HH^4(A)\times HH^4(A)\stackrel{0}{\rightarrow}HH^8(A)$} This follows by degree argument: $\deg HH^4(A)\geq -h$, $\deg HH^8(A)=-2h-2<-2h$.
\subsection{$HH^4(A)\times HH^5(A)\stackrel{0}{\rightarrow}HH^9(A)$} This is clear for $Q=D_{n+1}$, $n$ odd, $Q=E_7,\,E_8$ where $HH^9(A)=K[-2h-2]=0$.

Let $Q=D_{n+1}$, $n$ even or $Q=E_6$. Let $a\in HH^4(A)$, $b\in HH^5(A)$.The product $HH^2(A)\times HH^3(A)\rightarrow HH^5(A)$, $(x,y)\mapsto \langle x,y\rangle \zeta_0$ induces a nondegenerate bilinear form $\langle-,-\rangle$. If $ab\in HH^9(A)=HH^3(A)[-2h]$ is nonzero, then we can find a $c\in HH^2(A)$, such that $c(ab)=\zeta_0$. But this equals $(ca)b=0$ since $HH^2(A)\times HH^4(A)\stackrel{0}{\rightarrow}$ which gives us a contradiction.

\section{$HH^5(A)\times HH^5(A)\rightarrow HH^{10}(A)$}
 \begin{proposition}
  The multiplication of the subspace $U[-2]^*$ with $HH^5(A)$ is zero.
  
  The pairing on $Y^*[-h-2]$ is
  \begin{equation}\label{Y*xY*}
  \begin{array}{rcl}
   Y^*[-h-2]\times Y^*[-h-2]&\rightarrow& HH^{10}(A),\\
   (a,b)&\mapsto&\Omega(a,b)\varphi_4(\zeta_0),
   \end{array}
  \end{equation} 
  where the skew-symmetric bilinear form $\Omega(-,-)$ is given by the matrix $-M_\beta$ from subsection \ref{H1xH5}.
 \end{proposition}
\begin{proof}
We have $\deg HH^5(A)\geq-h-2$ and $\deg HH^{10}(A)\leq -2h-4$, so we can get a nonzero multiplication only by pairing bottom degree parts of $HH^5(A)$ which is $Y^*[-h-2]$. The product lies in the top degree part of $HH^{10}(A)=HH^4(A)[-2h]$ which is spanned by $\varphi_4(\zeta_0$). This gives us the pairing of the form (\ref{Y*xY*}).

We want to find the matrix $(\Omega(\varepsilon_i,\varepsilon_j))_{i,j}$ where $\varepsilon_i$ are a basis of \\$Y^*[-h-2]$, given in the section about $HH^5(A)$. Recall that the multiplication  $HH^1(A)\times HH^5(A)\rightarrow HH^6(A)$ was given by a skew-symmetric matrix $((M_\beta)_{i,j})_{i,j\in F}$, so that $\theta_0\varepsilon_i=\sum\limits_{k\in F}(M_\beta)_{k,i}\varphi_0(\omega_k)$.\\

We multiply $\varepsilon_i\varepsilon_j=\Omega(\varepsilon_i,\varepsilon_j)\varphi_4(\zeta_0)$ with $\theta_0$ (see \ref{thetaxzeta}):
\begin{equation}\label{Y*xY*1}
 \theta_0(\varepsilon_i\varepsilon_j)=\Omega(\varepsilon_i,\varepsilon_j)\varphi_5(\psi_0).
\end{equation}
Using associativity, this equals
\begin{equation}\label{Y*xY*2}
 (\theta_0\varepsilon_i)\varepsilon_j=\sum\limits_{k\in F}(M_\beta)_{k,i}\varphi_0(\omega_k)\varepsilon_j=(M_\beta)_{j,i}\psi_0=-(M_\beta)_{i,j}\varphi_5(\psi_0).
\end{equation}
We see from equations (\ref{Y*xY*1}) and (\ref{Y*xY*2}) that
\[
 \Omega(\varepsilon_i,\varepsilon_j)=-(M_\beta)_{i,j}.
\]

\end{proof}

\end{document}